# CONDITIONAL LEAST SQUARES ESTIMATION IN NONSTATIONARY NONLINEAR STOCHASTIC REGRESSION MODELS

By Christine Jacob[1]

*National Agronomical Research Institute (INRA)*

Let $\{Z_n\}$ be a real nonstationary stochastic process such that $E(Z_n|\mathcal{F}_{n-1}) \overset{\text{a.s.}}{<} \infty$ and $E(Z_n^2|\mathcal{F}_{n-1}) \overset{\text{a.s.}}{<} \infty$, where $\{\mathcal{F}_n\}$ is an increasing sequence of $\sigma$-algebras. Assuming that $E(Z_n|\mathcal{F}_{n-1}) = g_n(\theta_0, \nu_0) = g_n^{(1)}(\theta_0) + g_n^{(2)}(\theta_0, \nu_0)$, $\theta_0 \in \mathbb{R}^p$, $p < \infty$, $\nu_0 \in \mathbb{R}^q$ and $q \leq \infty$, we study the asymptotic properties of $\widehat{\theta}_n := \arg\min_\theta \sum_{k=1}^n (Z_k - g_k(\theta, \widehat{\nu}))^2 \lambda_k^{-1}$, where $\lambda_k$ is $\mathcal{F}_{k-1}$-measurable, $\widehat{\nu} = \{\widehat{\nu}_k\}$ is a sequence of estimations of $\nu_0$, $g_n(\theta, \widehat{\nu})$ is Lipschitz in $\theta$ and $g_n^{(2)}(\theta_0, \widehat{\nu}) - g_n^{(2)}(\theta, \widehat{\nu})$ is asymptotically negligible relative to $g_n^{(1)}(\theta_0) - g_n^{(1)}(\theta)$. We first generalize to this nonlinear stochastic model the necessary and sufficient condition obtained for the strong consistency of $\{\widehat{\theta}_n\}$ in the linear model. For that, we prove a strong law of large numbers for a class of submartingales. Again using this strong law, we derive the general conditions leading to the asymptotic distribution of $\widehat{\theta}_n$. We illustrate the theoretical results with examples of branching processes, and extension to quasi-likelihood estimators is also considered.

**1. Introduction.** Let $\{Z_n\}_{n \in \mathbb{N}}$ be an observed one-dimensional real stochastic process defined on a probability space $(\Omega, \mathcal{F}, P_{\theta_0, \nu_0})$ dependent on a unknown parameter $(\theta_0, \nu_0)$, $\theta_0 \in \Theta \subset \mathbb{R}^p$, $0 < p < \infty$, $\nu_0 \in \mathcal{N} \subset \mathbb{R}^q$, $0 \leq q \leq \infty$, and assumed to satisfy

$$M_Z: \forall n \quad E_{\theta_0, \nu_0}(Z_n|\mathcal{F}_{n-1}) = g_n(\theta_0, \nu_0) = g_n^{(1)}(\theta_0) + g_n^{(2)}(\theta_0, \nu_0),$$

$$E_{\theta_0, \nu_0}(Z_n^2|\mathcal{F}_{n-1}) \overset{\text{a.s.}}{<} \infty,$$

Received June 2008; revised July 2009.
[1]Supported in part by the French cooperation network ECONET.

where $\{\mathcal{F}_n\}$ is an increasing sequence of $\sigma$-algebras depending only on observed processes, $\theta_0$ is a unknown parameter that we want to estimate, $\nu_0$ is a nuisance parameter defined, when $q = \infty$, by $\nu_0 = \{\nu_{0n}\}$ with $g_n^{(2)}(\theta_0, \nu_0) = g_n^{(2)}(\theta_0, \nu_{0n})$, $g_n^{(1)}(\theta_0)$ is the $\mathcal{F}_{n-1}$-measurable parametric part of the model that may be nonlinear in $\theta_0$, and $g_n^{(2)}(\theta_0, \nu_0)$ is $\mathcal{F}_{n-1}$-measurable and such that $g_n^{(2)}(\theta_0, \nu_0) - g_n^{(2)}(\theta, \nu_0)$ is asymptotically negligible relative to $g_n^{(1)}(\theta_0) - g_n^{(1)}(\theta)$ (definition given in Section 4). The simplest example of asymptotic negligibility is when $g_n^{(2)}(\theta_0, \nu_0)$ is independent of $\theta_0$. The case $q = 0$ is defined by $g_n^{(2)}(\theta_0, \nu_0) = 0$, for all $n$, and corresponds to the classical parametric setting.

Examples of models $M_Z$ are nonlinear regression models with random covariates and heteroscedastic variances, stochastic dynamical models in discrete time, nonlinear time series model (TARMA, SETAR, bilinear processes), financial models (ARCH, GARCH and others) and branching processes, provided that the first two conditional moments at each $n$ of all of these processes are finite. This means, in particular, that processes with a heavy-tailed distribution (see [10] for such an example) do not belong to this class. However, a solution when $E_{\theta_0,\nu_0}(Z_n|\mathcal{F}_{n-1}) \stackrel{\text{a.s.}}{<} \infty$ with $\sigma_n^2(\theta_0, \nu_0) \stackrel{\text{a.s.}}{=} \infty$, where $\sigma_n^2(\theta_0, \nu_0) := E_{\theta_0,\nu_0}([Z_n - E_{\theta_0,\nu_0}(Z_n|\mathcal{F}_{n-1})]^2|\mathcal{F}_{n-1})$, could be to deal with the truncated process $\widetilde{Z}_n := Z_n 1_{\{Z_n \in I_n\}}$, where $\lim_n I_n \stackrel{\text{a.s.}}{=} \mathbb{R}$ since by defining $\widetilde{\eta}_n := \widetilde{Z}_n - E_{\theta_0,\nu_0}(\widetilde{Z}_n|\mathcal{F}_{n-1})$ and $\widetilde{g}_n(\theta_0, \nu_0) := E_{\theta_0,\nu_0}(\widetilde{Z}_n|\mathcal{F}_{n-1})$, we can then define $\widetilde{g}_n^{(1)}(\theta_0) = E_{\theta_0,\nu_0}(Z_n|\mathcal{F}_{n-1})$ and $\widetilde{g}_n^{(2)}(\theta_0, \nu_0) = -E_{\theta_0,\nu_0} \times (Z_n 1_{\{Z_n \notin I_n\}}|\mathcal{F}_{n-1})$.

We consider the class of weighted CLSE (conditional least squares estimators) of $\theta_0$ in the approximate model $\{g_k(\theta_0, \widehat{\nu})\}_{k \leq n}$, where $\widehat{\nu} = \{\widehat{\nu}_n\}$, $\widehat{\nu}_n$ being any estimation of $\nu_0$ based on observations up to $n$. We will consider two different settings:

$$\text{A1:} \quad \forall n \quad g_n(\theta, \widehat{\nu}) = g_n(\theta, \widehat{\nu}_n) \quad \text{or} \quad g_n(\theta, \widehat{\nu}) = g_n(\theta, \widehat{\nu}_{n_0});$$

$$\text{A2:} \ \forall n, \forall k \leq n \quad g_k(\theta, \widehat{\nu}) = g_k(\theta, \widehat{\nu}_n).$$

Such an estimator is defined by

$$(1.1) \quad \widehat{\theta}_n := \arg\min_{\theta \in \Theta} S_{n|\widehat{\nu}}(\theta), \qquad S_{n|\widehat{\nu}}(\theta) := \sum_{k=1}^{n} (Z_k - g_k(\theta, \widehat{\nu}))^2 \lambda_k^{-1},$$

where $\lambda_k$ is an $\mathcal{F}_{k-1}$-measurable variable independent of $(\theta_0, \nu_0)$. When $g_k(\theta, \widehat{\nu})$ has a first derivative $\mathbf{g}'_k(\theta, \widehat{\nu})$ in $\theta$, (1.1) implies that $\widehat{\theta}_n$ is an estimating equations estimator (EEE), that is,

$$(1.2) \quad \mathbf{Q}_{n|\widehat{\nu}}(\widehat{\theta}_n; \{Z_k, \mathbf{a}_k(\theta)\}) = 0^{p \times 1}, \qquad \{\mathbf{a}_k(\theta) \mathcal{F}_{k-1}\text{-measurable}\},$$

$$(1.3) \quad \mathbf{Q}_{n|\widehat{\nu}}(\theta; \{Z_k, \mathbf{a}_k(\theta)\}) := -\sum_{k=1}^{n}(Z_k - g_k(\theta, \widehat{\nu}))\mathbf{a}_k(\theta),$$

where, here, $\mathbf{a}_k(\theta) = \mathbf{g}'_k(\theta, \widehat{\nu})\lambda_k^{-1}$ for all $k$.

In the classical parametric setting $q = 0$, the weighted CLSE's and, more generally, the EEE's, are well studied and are known to have interesting properties. These estimators are robust to the form of distribution of the respective residuals $\{Z_n - g_n(\theta_0, \nu_0)\}$ since they require at most the knowledge of the first two conditional moments of the process at each time $n$ and their computation may be achieved, even in the case of complex or unknown likelihoods. When $\sigma_n^2(\theta_0, \nu_0)$ is an explicit function of $(\theta_0, \nu_0)$, such an estimator may be used to derive the empirical distribution of the estimated residuals $\{[Z_k - g_k(\widehat{\theta}_n, \widehat{\nu})][\sigma_k(\widehat{\theta}_n, \widehat{\nu})]^{-1}\}$, thereby allowing a re-estimation of $\theta_0$ by a maximum likelihood estimator (MLE) when the theoretical distribution may be modeled by a function of $\theta_0$, provided that the estimator is close enough to $\theta_0$ and is therefore strongly consistent [26]. In the particular setting $\sigma_n^2(\theta_0, \nu_0) = \sigma^2(\theta_0, \nu_0)\lambda_n$, where $\lambda_n$ is $\mathcal{F}_{n-1}$-measurable and independent of $(\theta_0, \nu_0)$, the optimal CLSE of $\theta_0$, from the point-of-view of the asymptotic variance, is obtained by making the errors of the model stationary, that is, by minimizing $\sum_{k=1}^{n}(Z_k - g_k(\theta, \nu_0))^2\lambda_k^{-1}$, and is equal to the optimal EEE, called the quasi-likelihood estimator (QLE) (see [15] for the optimality of the convergence rate in the branching process setting and [7] for the QLE). In the general case, if $g_n(\theta, \nu_0)$ has a first derivative in $\theta$, then a possible estimator is obtained by replacing $\sigma_k^2(\theta_0, \nu_0)$ by $\sigma_k^2(\theta, \widehat{\nu})$ in $\mathbf{Q}_{n|\widehat{\nu}}(\theta; \{Z_k\sigma_k^{-1}(\theta_0, \nu_0), \mathbf{g}'_k(\theta, \widehat{\nu})\sigma_k^{-1}(\theta_0, \nu_0)\})$. When $q = 0$, the obtained estimator is the QLE. This estimator is optimal from its asymptotic variance point-of-view within the class of estimators which solve (1.2) and (1.3), [7], and is moreover equal to the maximum likelihood estimator (MLE) when the conditional distribution of $Z_n$ belongs to an exponential family at each $n$ [32]. Another possible estimator is the weighted CLSE defined by (1.1), where $\{\lambda_k\}$ is a sequence of $\mathcal{F}_{k-1}$-measurable estimators of $\{\sigma_k^2(\theta_0, \nu_0)\}$ up to a multiplicative constant. Since

$$\min_{k \leq n} \frac{\sigma_k^2(\theta_0, \nu_0)}{\lambda_k} < \frac{\sum_{k=1}^{n}(Z_k - g_k(\theta, \widehat{\nu}))^2\lambda_k^{-1}}{\sum_{k=1}^{n}(Z_k - g_k(\theta, \widehat{\nu}))^2\sigma_k^{-2}(\theta_0, \nu_0)} < \max_{k \leq n} \frac{\sigma_k^2(\theta_0, \nu_0)}{\lambda_k},$$

if $\{\lambda_n\}$ is such that

$$(1.4) \quad 0 \stackrel{\text{a.s.}}{<} \varliminf_{n} \sigma_n^2(\theta_0, \nu_0)\lambda_n^{-1} \leq \varlimsup_{n} \sigma_n^2(\theta_0, \nu_0)\lambda_n^{-1} \stackrel{\text{a.s.}}{<} \infty,$$

then the asymptotic behavior of $\sum_{k=1}^{n}(Z_k - g_k(\theta, \widehat{\nu}))^2\lambda_k^{-1}$, and therefore of its argmin, should be close to that of $\sum_{k=1}^{n}(Z_k - g_k(\theta, \widehat{\nu}))^2\sigma_k^{-2}(\theta_0, \nu_0)$ and of its arg min. Finally, if $\widehat{F}_{n-1}$ denotes the set of $\mathcal{F}_{n-1}$-measurable variables, since $E_{\theta_0,\nu_0}(Z_n|\mathcal{F}_{n-1})$ is the best predictor of $Z_n$ based on $\widehat{F}_{n-1}$ in the least squares sense because $E_{\theta_0,\nu_0}(Z_n|\mathcal{F}_{n-1}) = \arg\min_{g \in \widehat{F}_{n-1}} E_{\theta_0,\nu_0}((Z_n - g)^2 \lambda_n^{-1}|\mathcal{F}_{n-1})$, provided that $\lambda_n$ is an $\mathcal{F}_{n-1}$-measurable variable independent of $g$ [9], a weighted CLSE should easily be strongly consistent. Since we are particularly interested in such a property, which is necessary when accurate knowledge of the true parameter is required, we will focus here on the asymptotic properties (strong consistency, asymptotic distribution), as $n \to \infty$, of the weighted CLSE solution of (1.1) in the general setting $M_Z$ with the weakest possible assumptions on the process behavior, conditionally on $\{\widehat{\nu}\}$. However, the results could easily be generalized to the QLE when a primitive of the estimating functions exists (see Section 9).

From now on, to simplify notation when studying $\{\widehat{\theta}_n\}$ which solve (1.1), we will use the normalized process $Y_n := Z_n \lambda_n^{-1/2}$ and denote by $M_Y$ the assumptions concerning $\{Y_n\}$ when $\{Z_n\}$ verifies $M_Z$. More precisely, defining $f_n^{(1)}(\theta_0) := g_n^{(1)}(\theta_0)\lambda_n^{-1/2}$ and $f_n^{(2)}(\theta_0, \nu_0) := g_n^{(2)}(\theta_0, \nu_0)\lambda_n^{-1/2}$, we have

$$M_Y\colon \forall n \quad E_{\theta_0,\nu_0}(Y_n|\mathcal{F}_{n-1}) = f_n(\theta_0, \nu_0) = f_n^{(1)}(\theta_0) + f_n^{(2)}(\theta_0, \nu_0),$$

$$E_{\theta_0,\nu_0}(Y_n^2|\mathcal{F}_{n-1}) \stackrel{\text{a.s.}}{<} \infty$$

and $\eta_n := Y_n - E_{\theta_0,\nu_0}(Y_n|\mathcal{F}_{n-1})$ is a martingale difference. Equations (1.1) and (1.4) are now written, respectively,

$$(1.5) \qquad \widehat{\theta}_n = \arg\min_{\theta \in \Theta} S_{n|\widehat{\nu}}(\theta), \qquad S_{n|\widehat{\nu}}(\theta) := \sum_{k=1}^{n}(Y_k - f_k(\theta, \widehat{\nu}))^2,$$

$$(1.6) \qquad 0 < \varliminf_n \sigma_n^2 \stackrel{\text{a.s.}}{\leq} \varlimsup_n \sigma_n^2 \stackrel{\text{a.s.}}{<} \infty, \qquad \sigma_n^2 := E(\eta_n^2|\mathcal{F}_{n-1}).$$

Among published works on the estimator consistency in $M_Y$, only the case $q = 0$ is considered and two large classes of proofs exist. One class is based on the stationarity and ergodicity assumptions of the process [26], on the strongest assumption of independence of the errors (classical regression) or on the explicit expression of the estimator according to the process together with the knowledge of its asymptotic behavior. It is, in particular, the case of a branching process when the corresponding model is linear in $\theta_0$ [8, 35]. The other approach is based on the (much more general) martingale difference property of $\{\eta_n\}$. Here, we are interested in this second class, which is particularly useful for processes. When dealing with the parametric linear model $f_n(\theta_0) = \theta_0^T \mathbf{W}_n$, where $\mathbf{W}_n$ is either a deterministic vector and $\{\eta_n\}$ are i.i.d. [20] or $\mathbf{W}_n$ is stochastic with $p = 1$ [21], then $\lim_n \widehat{\theta}_n \stackrel{\text{a.s.}}{=} \theta_0$ if

and only if

$$\lim_n \lambda_{\min}\left\{\sum_{k=1}^n \mathbf{W}_k \mathbf{W}_k^T\right\} \stackrel{\text{a.s.}}{=} \infty, \tag{1.7}$$

$\lambda_{\min}\{\sum_{k=1}^n \mathbf{W}_k \mathbf{W}_k^T\}$ being the smallest eigenvalue of $\sum_{k=1}^n \mathbf{W}_k \mathbf{W}_k^T$. Defining

$$D_n(\theta) := \sum_{k=1}^n [d_k(\theta)]^2, \qquad d_k(\theta) := f_k(\theta_0) - f_k(\theta), \tag{1.8}$$

(1.7) is equivalent to $\lim_n D_n(\theta) \stackrel{\text{a.s.}}{=} \infty$ for all $\theta \neq \theta_0$. This quantity is the identifiability criterion of $\theta_0$ in the model. It is interesting to observe that $\{\widehat{\theta}_n\}$ cannot be consistent, or even weakly consistent, on the set $\{\lim_n D_n(\theta) \stackrel{\text{a.s.}}{<} \infty\}$.

However, in the general nonlinear stochastic setting $M_Y$ with $q = 0$, under some Lipschitz property of the model, all published theorems of consistency require, besides the condition $\overline{\lim}_n \sigma_n^2 \stackrel{\text{a.s.}}{<} \infty$ and a condition of the type $\lim_n D_n(\widetilde{\theta}_n) \stackrel{\text{a.s.}}{=} \infty$ for some sequence $\{\widetilde{\theta}_n\} \in \Theta \setminus \theta_0$ (depending on the author), additional conditions concerning some rate of convergence to $\infty$ of $\{D_n(\cdot)\}$. Moreover, these conditions differ from one author to another ([1, 16, 18, 21–23, 31, 33, 36]; see [31] or [14] for some examples of models that do not verify these additional conditions).

Here, we generalize the necessary and sufficient condition (1.7) to our general nonstationary nonlinear stochastic model $M_Y$ with $0 \leq q \leq \infty$. When $q = 0$, we prove the strong consistency of $\{\widehat{\theta}_n\}$ on the set

$$\text{LIP}_\theta(\{f_k(\theta)\}) \cap \text{SI}_\theta(\{D_n(\theta)\}) \cap \text{VAR}_\theta(\{\sigma_k^2, d_k(\theta), D_k(\theta)\}), \tag{1.9}$$

where, in the following, "$\forall \delta > 0$" means "$\forall \delta > 0$ small enough" and

- $\text{LIP}_\theta(\{f_k(\theta)\})$ is the set of trajectories satisfying the following Lipschitz condition: for all $k$, there exists a nonnegative $\mathcal{F}_{k-1}$-measurable function $g_k$ and a function $h(\cdot): \mathbb{R}^+ \to \mathbb{R}^+$ with $\lim_{x \searrow 0} h(x) = 0$ such that for all $\theta_1 \in \Theta, \theta_2 \in \Theta, |f_k(\theta_1) - f_k(\theta_2)| \stackrel{\text{a.s.}}{\leq} h(\|\theta_1 - \theta_2\|) g_k$, where $\|\cdot\|$ is any norm in $\mathbb{R}^p$;
- $\text{VAR}_\theta(\{\sigma_k^2, d_k(\theta), D_k(\theta)\})$ that generalizes $\overline{\lim}_n \sigma_n^2 \stackrel{\text{a.s.}}{<} \infty$ (proved in Section 5) is the set

$$\left\{\forall \delta > 0, \sup_{\|\theta - \theta_0\| \geq \delta} \sum_{k=1}^\infty \sigma_k^2 [d_k(\theta)]^2 [D_k(\theta)]^{-2} \stackrel{\text{a.s.}}{<} \infty\right\}; \tag{1.10}$$

- $SI_\theta(\{D_n(\theta)\})$ concerns the identifiability condition generalizing (1.7):

$$\forall \delta > 0 \quad \inf_{\|\theta-\theta_0\|\geq \delta} D_n(\theta) \text{ is } \mathcal{F}_{n-1}\text{-measurable},$$

$$\lim_n \inf_{\|\theta-\theta_0\|\geq \delta} D_n(\theta) \stackrel{\text{a.s.}}{=} \infty.$$

The same terms $\text{LIP}_\theta(\cdot)$, $\text{VAR}_\theta(\cdot)$, $SI_\theta(\cdot)$ will indicate both the set of trajectories and the corresponding conditions verified by these sets. The result is then generalized to the setting $0 \leq q \leq \infty$, replacing, in each condition of (1.9), $f_k(\theta)$ by $f_k^{(1)}(\theta)$. This consistency result is due to an original SLLNSM (strong law of large numbers for submartingales). In addition, we show that the asymptotic distribution of the CLSE is easily derived from a classical CLT (central limit theorem), thanks to this SLLNSM. Therefore, this SLLNSM is the key result of this work.

The paper is organized in the following way. In Section 2, we give some examples of processes $\{Z_n\}$ satisfying $M_Z$. We deal with the consistency of $\{\widehat{\theta}_n\}$ in Section 3 when $q = 0$, and in Section 4 in the more general setting $0 \leq q \leq \infty$. This result is obtained thanks to an SLLNSM that is proved in Section 5 using submartingale properties [9], analytical lemmas and Wu's lemma concerning the consistency of estimators minimizing a contrast [34].

The consistency result obtained in the general setting $0 \leq q \leq \infty$ shows the robustness of this property with respect to the chosen model since, if $\{\widehat{\theta}_n\}$ is strongly consistent in a given model, then it is strongly consistent in every model "close," from the identifiability point-of-view, to this given model.

In Section 6, we give general conditions for obtaining the asymptotic distribution of $\{\widehat{\theta}_n\}$ from the classical CLT for martingales or for random sums. As in the classical nonlinear deterministic regression model [34], the proof is based on the Taylor series expansion of $\partial S_{n|\widehat{\nu}}(\theta)/\partial \theta$, where the convergence to 0 of the remaining term of the Taylor series is a direct consequence of the SLLNSM.

In Section 7, we estimate the part of $\nu_0$ involved in the asymptotic distribution of $\{\widehat{\theta}_n\}$ and give conditions for its consistency.

In Section 8, we give some examples in the single-type branching processes field. These processes model population dynamics. The population size $N_n$ at $n$ is defined by $N_n = \sum_{i=1}^{N_{n-1}} X_{n,i}$, where the offspring sizes $\{X_{n,i}\}_i$, given $\mathcal{F}_{n-1}$, are i.i.d. with mean $m_{\theta_0,\nu_0}(F_{n-1})$ and variance $\sigma^2_{\theta_0,\nu_0}(F_{n-1})$, $F_{n-1}$ denoting the set of random variables involved in $\mathcal{F}_{n-1}$, and $\mathcal{F}_{n-1}$ being generated by $\{N_k\}_{k\leq n-1}$ and possibly environmental processes until $n$. Estimation in this field is well understood in the framework of a BGW (Bienaymé–Galton–Watson) process $m_{\theta_0,\nu_0}(F_{n-1}) = m_0$, $\sigma^2_{\theta_0,\nu_0}(F_{n-1}) = \sigma_0^2$ or a derived process (BGW with immigration, controlled branching processes), assuming a linear model in $\theta_0$, and are predominantly based on the observation of the size $N_n$ of the population at each time $n$. The estimators are most often moment estimators, MLE, CLSE, QLE. The asymptotic properties as $n \to \infty$ are derived on the nonextinction set from the explicit expression for the estimator according to $\{N_n\}$ using the asymptotic behavior of the process, when suitably normalized. An overview of the references may be found in [8, 35], or in [25] for additional references for multitype processes.

Much more difficult is the study in the nonlinear case when there is no explicit expression for the estimators. The MLE is generally not easily computed because of large combinatorial terms. However, if we write $Y_n := N_n N_{n-1}^{-1/2} = m_{\theta_0, \nu_0}(F_{n-1}) N_{n-1}^{1/2} + \eta_n$, where $\eta_n = N_{n-1}^{-1/2} \sum_{i=1}^{N_{n-1}} (X_{n,i} - m_{\theta_0, \nu_0}(F_{n-1}))$, then $\eta_n$ is approximately normally distributed on the nonextinction set as soon as $n$ is large enough. Consequently, the CLSE with $\{\lambda_n\}$ satisfying (1.4) is approximately equal to the MLE on this set.

Therefore, we began to study the CLSE of $\theta_0$ in size-dependent branching processes [15, 24] and in regenerative branching processes [16]. The results presented in this paper improve on, and generalize, these results. The first example is a supercritical single-type BGW process, where the offspring mean $m_0$ is estimated. We give the asymptotic properties of $\widehat{m}_n$. The results are well known [3, 8], but the indirect form of proof given here for the consistency does not require accurate knowledge of the asymptotic behavior of the process, as is the case in the classical direct proof based on the analytical expression of the estimator. We then deal with some near-critical size-dependent branching processes which model the amplification process in the polymerase chain reaction setting [27]. In this setting, we prove the strong consistency and asymptotic distribution of the CLSE of the parameters of the replication probability for each model. The more complex model contains an explosive part tending to $\infty$, a persistent bounded part and a transient part. In this model, despite these three different kinds of behavior, we prove the consistency and the asymptotic normality of the CLSE of a bidimensional parameter with components belonging to the persistent part and the transient one.

Finally, in Section 9, we extend the consistency conditions to estimators solving (1.2), (1.3).

In the following sections, we will simply write $E(\cdot)$ for $E_{\theta_0, \nu_0}(\cdot)$.

**2. Examples of models.** Here, we give some examples of classical models satisfying $M_Z$.

2.1. *Regression models.* The observed variable $Z_n$ is explained by a parametric regression function $g_n(\theta_0)$ of a random (or deterministic) vector of observed covariates or coprocesses $\{\mathbf{X}_k\}_{k \leq n}$ and the residuals $\{Z_n - g_n(\theta_0)\}$ are assumed independent with finite variances. The time regression model with $d = 1$ and $X_{n,l} = n$ is a simple example.

The general class $M_Z$ allows the extension of these models to $g_n(\theta_0, \nu_0)$ and nonindependent residuals such as the following example: $Z_n = g_n(\theta_0) + \epsilon_{n-1}\epsilon_n$. Since we may write $\epsilon_{n-1}$ as a function of $\{\{Z_k - g_k(\theta_0)\}_{k=1}^{n-1}, \epsilon_0\}$, this model belongs to the $M_Z$ class.

2.2. *Financial time series.* The most well known of these are the ARCH model introduced by Engle [5] and the more general GARCH models. Let $\xi_n = s_n(\theta_0)U_n$, where the $\{U_n\}$ are i.i.d. $(0,1)$, $s_n(\theta_0) \geq 0$, $\mathcal{F}_{n-1}$ is generated by $\{\xi_k^2\}_{k \leq n-1}$ and $s_n^2(\theta_0)$ is $\mathcal{F}_{n-1}$-measurable with $s_n^2(\theta) = \alpha_0 + \sum_{j=1}^q \alpha_j \xi_{n-j}^2 + \sum_{j=1}^p \beta_j s_{n-j}^2(\theta)$. The process $\{s_n(\theta_0)\}$ is called *volatility*. Then $E(\xi_n^2|\mathcal{F}_{n-1}) = s_n^2(\theta_0)$ and $\{\xi_n\}$ follows a GARCH$(p,q)$ model. If $\{\xi_n\}$ is observed, then $Z_n$ is defined in the following way: $Z_n = \xi_n^2 = s_n^2(\theta_0) + s_n^2(\theta_0)(U_n^2 - 1)$, implying that $g_n(\theta_0) = s_n^2(\theta_0)$, $\sigma_n(\theta_0, \nu_0) = s_n^2(\theta_0)[E((U_n^2-1)^2|\mathcal{F}_{n-1})]^{1/2}$. More generally, $\{\xi_n\}$ may be nondirectly observed [14, 26].

2.3. *Linear or nonlinear time series models that may depend on nonstationary exogenous inputs $\{u_n\}$, such as ARMAX models.* We define $e_n := Z_n - g_n(\theta_0)$, where $g_n(\theta) = \sum_{k=1}^p \alpha_k Z_{n-k} + \sum_{k=1}^q \beta_k e_{n-k} + \sum_{k=0}^b \gamma_k u_{n-k}$. Here, $\mathcal{F}_{n-1}$ is generated by $\{Z_k\}_{k \leq n-1}, \{u_k\}_{k \leq n}$.

2.4. *Single-type discrete-time branching processes.* This class of models is described in Sections 1 and 8.

2.5. *Models with observation errors.* In practice, the model may be $Z_n^{\text{th}} = g_n^{\text{th}}(\theta_0) + e_n$, $Z_n = Z_n^{\text{th}} + u_n$, where $Z_n$ is the observation of the theoretical unknown variable $Z_n^{\text{th}}$, $u_n$ is the observation error with a unknown distribution and $g_n^{\text{th}}(\theta_0) = E(Z_n^{\text{th}}|\mathcal{F}_{n-1}^{\text{th}})$. Then, assuming that $E(u_n|\mathcal{F}_{n-1}) = 0$, we have $g_n(\theta_0, \nu_0) = E(g_n^{\text{th}}(\theta_0)|\mathcal{F}_{n-1})$. For example, consider the BGW process $\{N_n^{\text{th}}\}$ with $N_n = N_n^{\text{th}} + u_n$, $f_n^{\text{th}}(m) = m(N_{n-1}^{\text{th}})^{1/2}$. Then, using the first order Taylor series expansion, we get

$$f_n(m_0, \nu_0) := E(f_n^{\text{th}}(\theta_0)|\mathcal{F}_{n-1})$$
$$= mE([N_{n-1} - u_{n-1}]^{1/2}|\mathcal{F}_{n-1})$$
$$= mN_{n-1}^{1/2} + mO(E(u_{n-1}^2|\mathcal{F}_{n-1})N_{n-1}^{-3/2}),$$

where $\nu = \{mE(u_{n-1}^2|\mathcal{F}_{n-1})\}$, $f_n^{(2)}(m, \nu) = mO(E(u_{n-1}^2|\mathcal{F}_{n-1})N_{n-1}^{-3/2})$. The asymptotic negligibility of $\{f_n^{(2)}(m, \widehat{\nu})\}$ is obtained on the nonextinction set for any bounded sequence $\widehat{\nu}$ since $\lim_n N_n \stackrel{\text{a.s.}}{=} \infty$ on this set.

2.6. *Multivariate stochastic regression models.* Estimation in this field may also be expressed as a finite fixed set of one-dimensional models belonging to one of the previous types. Let us assume that $\mathbf{Z}_k \in \mathbb{R}^d$ with $E(\mathbf{Z}_n|\mathcal{F}_{n-1}) = \mathbf{g}_n(\theta_0, \nu_0)$ and let $\mathbf{\Sigma}_n$ be a known $\mathcal{F}_{n-1}$-measurable positive definite estimation of the conditional variance-covariance matrix of $\mathbf{Z}_n$ up to some unknown multiplicative constant. Then $\widehat{\theta}_n = \arg\min_{\theta \in \Theta} S_{n|\widehat{\nu}}(\theta)$, where

$$S_{n|\widehat{\nu}}(\theta) := \sum_{k=1}^{n} (\mathbf{Z}_k - \mathbf{g}_k(\theta, \widehat{\nu}))^T \mathbf{\Sigma}_k^{-1} (\mathbf{Z}_k - \mathbf{g}_k(\theta, \widehat{\nu})).$$

Since $\mathbf{\Sigma}_k$ is positive definite, we may write $\mathbf{\Sigma}_k = \mathbf{U}_k \mathbf{\Lambda}_k \mathbf{U}_k^{-1}$, where $\mathbf{U}_k$ is an orthogonal matrix and $\mathbf{\Lambda}_k$ is the diagonal matrix of the eigenvalues of $\mathbf{\Sigma}_k$. Therefore, writing $\mathbf{Y}_k = \mathbf{\Lambda}_k^{-1/2} \mathbf{U}_k^{-1} \mathbf{Z}_k$, where $\mathbf{Y}_k = (Y_{k,1}, \ldots, Y_{k,d})^T$, we get $\mathbf{f}_k(\theta_0, \widehat{\nu}) = \mathbf{\Lambda}_k^{-1/2} \mathbf{U}_k^{-1} \mathbf{g}_k(\theta_0, \widehat{\nu})$, $\mathbf{Y}_k - E(\mathbf{Y}_k|\mathcal{F}_{k-1})$ is a martingale difference and

(2.1) $$S_{n|\widehat{\nu}}(\theta) = \sum_{j=1}^{d} \sum_{k=1}^{n} (Y_{k,j} - f_{k,j}(\theta, \widehat{\nu}))^2.$$

**3. Consistency in $M_Y$ with $q = 0$.** In this section, we generalize (1.7) to the setting of the model $M_Y$ under $q = 0$. So, we have $d_k^{(1)}(\theta) = d_k(\theta)$. Let us assume that $\theta_0 \in \Theta$, where $\Theta$ is an open set and $\overline{\Theta}$ is compact. From now on (in all sections), let $\{\gamma \in B_\delta^c\} := \{\gamma : \|\gamma - \gamma_0\| \geq \delta\}$, where $\gamma$ may be any subset of $(\theta, \nu)$.

PROPOSITION 3.1. *Let $\Omega_\infty \subset \Omega$ be defined by*

$$\Omega_\infty = \mathrm{LIP}_\theta(\{f_k(\theta)\}) \cap \mathrm{SI}_\theta(\{D_n(\theta)\}) \cap \mathrm{VAR}_\theta(\{\sigma_k^2, d_k(\theta), D_k(\theta)\}).$$

*Let us assume that $P(\Omega_\infty) > 0$. Then $\lim_n \widehat{\theta}_n \stackrel{a.s.}{=} \theta_0$ on $\Omega_\infty$.*

REMARK 3.1. In the linear model, $f_n(\theta) = \theta^T \mathbf{W}_n$ with $\overline{\lim}_n \sigma_n^2 \stackrel{a.s.}{<} \infty$, $\Omega_\infty$ is reduced to $\{\lim_n \lambda_{\min}\{\sum_{k=1}^{n} \mathbf{W}_k \mathbf{W}_k^T\} \stackrel{a.s.}{=} \infty\}$ thanks to Proposition 5.2.

PROOF OF PROPOSITION 3.1. We use Wu's lemma ([34], see Lemma A.1) and Wu's decomposition based on $Y_k - f_k(\theta) = \eta_k + d_k(\theta)$ [34], implying $S_n(\theta) - S_n(\theta_0) = D_n(\theta) + 2L_n(\theta)$, where $L_n(\theta) = \sum_{k=1}^{n} \eta_k d_k(\theta)$ and, consequently,

(3.1) $$\inf_{\theta \in B_\delta^c} S_n(\theta) - S_n(\theta_0) \geq \inf_{\theta \in B_\delta^c} D_n(\theta) \Big[ 1 - 2 \sup_{\theta \in B_\delta^c} |L_n(\theta)| [D_n(\theta)]^{-1} \Big].$$

The proof then follows directly from the SLLNSM (Proposition 5.1) applied to $d_k(\theta) = f_k(\theta_0) - f_k(\theta)$ and $\widetilde{\Theta} = B_\delta^c$. □

Let us now assume that $\widehat{\theta}_{h,n} = \arg\min_{\theta \in \Theta} \sum_{k=h+1}^n (Y_k - f_k(\theta))^2$, where $h$ may depend on $n$ (e.g., $n - h$ is constant for all $n$) and let $L_n(\theta) - L_h(\theta) =: L_{h,n}(\theta)$ and $D_n(\theta) - D_h(\theta) =: D_{h,n}(\theta)$. Let us define

$$\mathrm{RAT}_\theta(D_n(\theta)[D_{h,n}(\theta)]^{-1}\}) := \Big\{\varlimsup_n \sup_{\theta \in B_\delta^c} D_n(\theta)[D_{h,n}(\theta)]^{-1} \overset{\mathrm{a.s.}}{<} \infty\Big\}.$$

PROPOSITION 3.2. *Let $\Omega_\infty \subset \Omega$ be defined by*

$$\Omega_\infty = \mathrm{LIP}_\theta(\{f_k(\theta)\}) \cap \mathrm{SI}_\theta(\{D_n(\theta)\}) \cap \mathrm{VAR}_\theta(\{\sigma_k^2, d_k(\theta), D_k(\theta)\})$$
$$\cap \mathrm{RAT}_\theta(\{D_n(\theta)[D_{h,n}(\theta)]^{-1}\}).$$

*Let us assume that $P(\Omega_\infty) > 0$. Then $\lim_n \widehat{\theta}_{h,n} \overset{\mathrm{a.s.}}{=} \theta_0$ on $\Omega_\infty$.*

PROOF. When $h$ is fixed, we are in the setting of Proposition 3.1. Therefore, let us assume that $h \to \infty$ as $n \to \infty$. According to (3.1) written with $S_{h,n}(\theta)$ instead of $S_n(\theta)$, it is sufficient to prove that $\lim_n \sup_\theta L_{h,n}(\theta)[D_{h,n}(\theta)]^{-1} = 0$ when $\lim_n \inf_{\theta \in B_\delta^c} D_{h,n}(\theta) \overset{\mathrm{a.s.}}{=} \infty$. For that, we use Proposition 5.1 with

$$\frac{L_{h,n}(\theta)}{D_{h,n}(\theta)} = \frac{L_n(\theta)}{D_n(\theta)} \frac{D_n(\theta)}{D_{h,n}(\theta)} - \frac{L_h(\theta)}{D_h(\theta)} \frac{D_h(\theta)}{D_{h,n}(\theta)}. \qquad \Box$$

Let us now assume, as in the last item of Section 2, the multidimensional case $\mathbf{Z}_n \in \mathbb{R}^d$. The CLSE of $\theta_0$ is then given by (2.1) and, more generally, by

$$\widehat{\theta}_{h,n} = \arg\min_{\theta \in \Theta} \sum_{j=1}^d \sum_{k=h+1}^n (Y_{k,j} - f_{k,j}(\theta))^2.$$

Let $D_{h,n,j}(\theta) = \sum_{k=h+1}^n (f_{k,j}(\theta) - f_{k,j}(\theta_0))^2$ and $D_{n,j}(\theta) := D_{0,n,j}(\theta)$.

PROPOSITION 3.3. *Let $\Omega_\infty \subset \Omega$ be defined by*

$$\Omega_\infty = \bigcap_j \{\mathrm{LIP}_\theta(\{f_{k,j}(\theta)\}) \cap \mathrm{SI}_\theta(\{D_{n,j}(\theta)\})$$
$$\cap \mathrm{VAR}_\theta(\{\sigma_{k,j}^2, d_{k,j}(\theta), D_{k,j}(\theta)\}) \cap \mathrm{RAT}_\theta(\{D_{n,j}(\theta)[D_{h,n,j}(\theta)]^{-1}\})\}.$$

*Let us assume that $P(\Omega_\infty) > 0$. Then $\lim_n \widehat{\theta}_{h,n} \overset{\mathrm{a.s.}}{=} \theta_0$ on $\Omega_\infty$.*

The proof follows directly from Proposition 3.2 applied to each $j = 1, \ldots, d$.

**4. Consistency in $M_Y$ containing a nuisance part.** Let us now assume that $\{f_n^{(2)}(\theta_0, \nu_0)\}$ is not identically null. We prove the consistency conditionally on a given sequence of estimations $\{\widehat{\nu}\}$.

As in Section 3, $\theta_0 \in \Theta$, where $\Theta$ is an open set and $\overline{\Theta}$ is compact, and, under A2, $\widehat{\nu}_n \in \mathcal{N}$, where $\overline{\mathcal{N}}$ is compact.

Let $D_n(\theta, \widehat{\nu}) := \sum_{k=1}^n [d_k(\theta, \widehat{\nu})]^2$, $D_n^{(i)}(\theta, \widehat{\nu}) := \sum_{k=1}^n [d_k^{(i)}(\theta, \widehat{\nu})]^2$, $d_k^{(i)}(\theta, \widehat{\nu}) := f_k^{(i)}(\theta_0, \nu_0) - f_k^{(i)}(\theta, \widehat{\nu})$, $L_n^{(i)}(\theta, \widehat{\nu}) := \sum_{k=1}^n \eta_k d_k^{(i)}(\theta, \widehat{\nu})$, $i = 1, 2$, and so on.

Let us define the following asymptotic negligibility property under A1:

$\mathrm{AN}_\theta(\{d_n^{(2)}(\theta, \widehat{\nu}), d_n^{(1)}(\theta)\})$: $\forall \delta > 0$

$$\varlimsup_n \left[\sup_{\theta \in B_\delta^c} |d_n^{(2)}(\theta, \widehat{\nu})|\right] \left[\inf_{\theta \in B_\delta^c} |d_n^{(1)}(\theta)|\right]^{-1} \stackrel{\text{a.s.}}{=} 0.$$

Under A2, we define $\mathrm{AN}_{\theta,\nu}(\{d_n^{(2)}(\theta, \nu), d_n^{(1)}(\theta)\})$ in the same way, replacing $\widehat{\nu}$ by $\nu$ and $\sup_{\theta \in B_\delta^c}$ by $\sup_{(\theta,\nu) \in B_\delta^c}$ in $\mathrm{AN}_\theta(\{d_n^{(2)}(\theta, \widehat{\nu}), d_n^{(1)}(\theta)\})$.

PROPOSITION 4.1. *Let $\Omega_\infty \subset \Omega$ be defined under* A1 *by*

$$\Omega_\infty = \mathrm{LIP}_\theta(\{f_k^{(1)}(\theta)\}) \cap \mathrm{LIP}_\theta(\{f_k^{(2)}(\theta, \widehat{\nu})\}) \cap \mathrm{SI}_\theta(\{D_n^{(1)}(\theta)\})$$

$$\cap \mathrm{AN}_\theta(\{d_n^{(2)}(\theta, \widehat{\nu}), d_n^{(1)}(\theta)\}) \cap \mathrm{VAR}_\theta(\{\sigma_k^2, d_k^{(1)}(\theta), D_k^{(1)}(\theta)\}).$$

*Under* A2, *replace $\widehat{\nu}$ by $\nu$ and $\mathrm{LIP}_\theta(\cdot)$, $\mathrm{AN}_\theta(\cdot)$ by $\mathrm{LIP}_{\theta,\nu}(\cdot)$, $\mathrm{AN}_{\theta,\nu}(\cdot)$, respectively. Let us assume that $P(\Omega_\infty) > 0$. Then $\lim_n \widehat{\theta}_n \stackrel{\text{a.s.}}{=} \theta_0$ on $\Omega_\infty$.*

REMARK 4.1. If, for all $k$,

$$\sup_{\theta_1, \theta_2 : f_k^{(1)}(\theta_1) \neq f_k^{(1)}(\theta_2)} |(f_k^{(2)}(\theta_1, \widehat{\nu}) - f_k^{(2)}(\theta_2, \widehat{\nu}))[f_k^{(1)}(\theta_1) - f_k^{(1)}(\theta_2)]^{-1}|$$

is $\mathcal{F}_{k-1}$-measurable, then $\mathrm{LIP}_\theta(\{f_k^{(2)}(\theta, \widehat{\nu})\})$ is satisfied under $\mathrm{LIP}_\theta(\{f_k^{(1)}(\theta)\})$.

PROOF OF PROPOSITION 4.1. We assume A1, implying that $f_k(\theta, \widehat{\nu})$ is independent of $n$ and is therefore $\mathcal{F}_{k-1}$-measurable. As in the proof of Proposition 3.1,

$$S_{n|\widehat{\nu}}(\theta) - S_{n|\widehat{\nu}}(\theta_0) = D_n(\theta, \widehat{\nu}) - D_n^{(2)}(\theta_0, \widehat{\nu}) + 2L_n(\theta, \widehat{\nu}) - 2L_n^{(2)}(\theta_0, \widehat{\nu}).$$

Since $D_n(\theta, \widehat{\nu}) = D_n^{(1)}(\theta) + D_n^{(2)}(\theta, \widehat{\nu}) + 2\sum_{k=1}^n d_k^{(1)}(\theta) d_k^{(2)}(\theta, \widehat{\nu})$, using Hölder's inequality, we get

$$\inf_{\theta \in B_\delta^c} S_{n|\widehat{\nu}}(\theta) - S_{n|\widehat{\nu}}(\theta_0)$$

$$
(4.1) \qquad \geq \inf_{\theta \in B_\delta^c} D_n^{(1)}(\theta) \left[ 1 - 2 \sup_{\theta \in B_\delta^c} \left[ \frac{D_n^{(2)}(\theta, \widehat{\nu})}{D_n^{(1)}(\theta)} \right]^{1/2} \right.
$$

$$
\left. - \sup_{\theta \in B_\delta^c} \frac{D_n^{(2)}(\theta_0, \widehat{\nu})}{D_n^{(1)}(\theta)} - 2 \sup_{\theta \in B_\delta^c} \frac{|L_n(\theta, \widehat{\nu}) - L_n^{(2)}(\theta_0, \widehat{\nu})|}{D_n^{(1)}(\theta)} \right].
$$

The result then follows from Wu's lemma A.1, the fact that $d_k(\theta, \widehat{\nu}) - d_k^{(2)}(\theta_0, \widehat{\nu}) = d_k^{(1)}(\theta) + d_k^{(2)}(\theta, \widehat{\nu}) - d_k^{(2)}(\theta_0, \widehat{\nu})$ and Proposition 5.1.

The proof under A2 is similar, replacing $\widehat{\nu}$ by $\nu$ and $\sup_\theta(\cdot)$ by $\sup_{\theta,\nu}(\cdot)$ in (4.1). □

Of course, Propositions 3.2 and 3.3 can also be easily generalized to $q > 0$.

## 5. Strong law of large numbers for submartingales.

PROPOSITION 5.1. *Let $\widetilde{\Theta} \subset \mathbb{R}^p$, $\overline{\widetilde{\Theta}}$ compact, $p < \infty$. Let $\{\mathcal{F}_k\}$ be an increasing sequence of $\sigma$-algebras on $\Omega$ and $L_n(\theta) = \sum_{k=1}^n \eta_k d_k(\theta)$, $\theta \in \widetilde{\Theta}$, where, for all $k$, $\eta_k$ is any $\mathcal{F}_k$-measurable variable such that $E(\eta_k|\mathcal{F}_{k-1}) = 0$, $E(\eta_k^2|\mathcal{F}_{k-1}) = \sigma_k^2 \stackrel{\text{a.s.}}{<} \infty$ and $d_k(\theta)$ is any $\mathcal{F}_{k-1}$-measurable variable. For all $k$, $n$, let $d_{*k}(\theta)$ be $\mathcal{F}_{k-1}$-measurable, $D_{*n}(\theta) = \sum_{k=1}^n d_{*k}^2(\theta)$, $D_n(\theta) = \sum_{k=1}^n d_k^2(\theta)$. Let $\Omega_\infty \subset \Omega$ be defined by*

$$
\Omega_\infty = \text{LIP}_\theta(\{d_k(\theta)\}) \cap \text{LIP}_\theta(\{d_{*k}(\theta)\}) \cap \text{SI}_\theta(\{D_{*n}(\theta)\})
$$
$$
\cap \text{VAR}_\theta(\{\sigma_k^2, d_k(\theta), D_{*k}(\theta)\}),
$$

*where $\text{LIP}_\theta(\cdot)$, $\text{SI}_\theta(\cdot)$ and $\text{VAR}_\theta(\cdot)$ are defined in Section 1. Let us assume that $P(\Omega_\infty) > 0$. Then*

$$
(5.1) \qquad \limsup_n \sup_{\theta \in \widetilde{\Theta}} |L_n(\theta)| [D_{*n}(\theta)]^{-1} \stackrel{\text{a.s.}}{=} 0 \qquad on\ \Omega_\infty.
$$

The proof is in the Appendix.

PROPOSITION 5.2. *Let $\{\eta_k, d_k(\theta), D_k(\theta)\}$ be as in Proposition 5.1. Let us assume that for each $\delta > 0$, there exists a random $k_{m,\delta} \in \mathbb{N}$ such that $k_{m,\delta} = \min\{k \geq 1 : [\inf_{\theta \in B_\delta^c} d_k(\theta)]^2 > 0\}$ exists. Then*

$$
\left\{ \overline{\lim_n} \sigma_n^2 \stackrel{\text{a.s.}}{<} \infty \right\} \subset \text{VAR}_\theta(\{\sigma_k^2, d_k(\theta), D_k(\theta)\}).
$$

PROOF. Let us assume that $\overline{\lim}_n \sigma_n^2 \stackrel{\text{a.s.}}{<} \infty$. For each $\theta$ and each trajectory, let us define $f(x) := [d_k(\theta)]^2$ for $x \in [(k-1), k[$. Then

$$
\sum_{k=k_{m,\delta}+1}^\infty \frac{d_k^2(\theta)}{[D_k(\theta)]^2} = \sum_{k=k_{m,\delta}+1}^\infty \frac{\int_{k-1}^k f(x)\,dx}{[\int_0^k f(u)\,du]^2}
$$

$$\leq \int_{k_{m,\delta}}^{\infty} \frac{f(x)\,dx}{[\int_0^x f(u)\,du]^2}$$

$$\leq \left[\frac{1}{\int_0^x f(u)\,du}\right]_\infty^{k_{m,\delta}}$$

$$\leq \frac{1}{\inf_{\theta \in B_\delta^c}[d_{k_{m,\delta}}(\theta)]^2}. \qquad \square$$

**6. Asymptotic distribution.** This section is devoted to general conditions leading to an asymptotic distribution of $\widehat{\theta}_n$ given $\widehat{\nu}$, where $\widehat{\theta}_n$ is solution of (1.5) in which either $f_n(\theta_0, \nu_0) = f_n^{(1)}(\theta_0) + f_n^{(2)}(\theta_0, \nu_0)$ or $f_n(\theta_0, \nu_0) = f_n^{(1)}(\theta_0, \nu_0)$. This latter case is suitable when $\{\widehat{\theta}_n, \widehat{\nu}_n\}$ is strongly consistent, but the assumptions leading to the asymptotic distribution are fulfilled only with respect to $\theta_0$ (see examples in Section 8.2).

We introduce the following notation [and similarly for $f_k(\theta, \widehat{\nu})$, $k \leq n$]:

$$\frac{\partial S_{n|\widehat{\nu}}}{\partial \theta_i}(\theta) =: S'_{n|\widehat{\nu};i}(\theta) =: \mathbf{S}'_{n|\widehat{\nu}}(\theta)[i],$$

$$\frac{\partial S'_{n|\widehat{\nu};i}}{\partial \theta_j}(\theta) =: S''_{n|\widehat{\nu};i,j}(\theta) =: \mathbf{S}''_{n|\widehat{\nu}}(\theta)[i,j],$$

where $\mathbf{S}'_{n|\widehat{\nu}}(\theta)$ [resp., $\mathbf{S}''_{n|\widehat{\nu}}(\theta)$] is a $p \times 1$ (resp., $p \times p$) matrix. Let $\mathbf{M}_n = [\sum_{k=1}^n \mathbf{f}'_k(\theta_0, \widehat{\nu})\mathbf{f}'^T_k(\theta_0, \widehat{\nu})]^{-1}$ (assumed to exist for all $n$ sufficiently large) and let $a_n \in (0,1)$, $\theta_n = \theta_0 + a_n(\widehat{\theta}_n - \theta_0)$. Let us define, for $q \leq \infty$, the following sets of trajectories:

$$\mathrm{UNC}(\{\theta_n\}): \left\{\lim_n \left[\sum_{k=1}^n \mathbf{f}'_k(\theta_n, \widehat{\nu})\mathbf{f}'^T_k(\theta_n, \widehat{\nu})\right]\mathbf{M}_n \stackrel{P}{=} \mathbf{I}\right\},$$

$$\mathrm{LIM}(\{\theta_n\}): \bigcap_{i,j,l}\left\{\lim_n \left|\sum_{k=1}^n (f_k(\theta_0, \nu_0) - f_k(\theta_n, \nu_0))f''_{k;i,l}(\theta_n, \widehat{\nu})\mathbf{M}_n[l,j]\right| \stackrel{P}{=} 0\right\}$$

$$\cap \left\{\lim_n \left|\sum_{k=1}^n (f_k(\theta_n, \nu_0) - f_k(\theta_n, \widehat{\nu}))f''_{k;i,l}(\theta_n, \widehat{\nu})\mathbf{M}_n[l,j]\right| \stackrel{P}{=} 0\right\}.$$

REMARK 6.1. Note that $\mathrm{UNC}(\{\theta_n\})$ may also be written as follows: for all $(i,j)$,

$$\lim_n \sum_l \sum_{h:\mathbb{N}_h(i,l) \neq \varnothing}\left(\sum_{k \leq n, k \in \mathbb{N}_h(i,l)} \frac{f'_{k;i}(\theta_n, \widehat{\nu})f'_{k;l}(\theta_n, \widehat{\nu})}{f'_{k;i}(\theta_0, \widehat{\nu})f'_{k;l}(\theta_0, \widehat{\nu})} f'_{k;i}(\theta_0, \widehat{\nu})f'_{k;l}(\theta_0, \widehat{\nu})\right)$$

$$\times \left( \sum_{k \le n, k \in \mathbb{N}_h(i,l)} f'_{k;i}(\theta_0, \widehat{\nu}) f'_{k;l}(\theta_0, \widehat{\nu}) \right)^{-1}$$

$$\times \sum_{k \le n, k \in \mathbb{N}_h(i,l)} f'_{k;i}(\theta_0, \widehat{\nu}) f'_{k;l}(\theta_0, \widehat{\nu}) \mathbf{M}_n[l, j] \stackrel{\text{a.s.}}{=} \delta_i(j),$$

where $\mathbb{N}_1(i, l) = \{k : f'_{k;i}(\theta_0, \widehat{\nu}) f'_{k;l}(\theta_0, \widehat{\nu}) > 0\}$, and similarly for $\mathbb{N}_2(i, l)$ with "<0" instead of ">0." Therefore, UNC($\{\theta_n\}$) is verified on the set

$$\bigcap_i \left\{ \lim_{n, \theta_n \to \theta_0} \sup_{k \,:\, f'_{k;i}(\theta_0, \widehat{\nu}) \ne 0} |f'_{k;i}(\theta_n, \widehat{\nu})[f'_{k;i}(\theta_0, \widehat{\nu})]^{-1} - 1| \stackrel{\text{a.s.}}{=} 0 \right\}.$$

PROPOSITION 6.1. *Let us assume that $f_k(\theta, \nu)$ has second derivatives in $\theta$ for each $k$ (and for $\nu = \widehat{\nu}$ under A1 and for each $\nu$ under A2) and that there exists a $p \times p$ $\mathcal{F}_{n-1}$-measurable matrix $\mathbf{\Psi}_n$ such that $P(\Omega_\infty) > 0$, where $\Omega_\infty$ is defined under A1 by*

$$\Omega_\infty = \bigcap_{i,j,l} \{ \mathrm{LIP}_\theta(\{f''_{k;i,l}(\theta, \widehat{\nu})\}) \cap \mathrm{SI}(\{(\mathbf{M}_n[l, j])^{-1}\})$$

$$\cap \mathrm{VAR}_\theta(\{\sigma_k^2, f''_{k,i,l}(\theta, \widehat{\nu}), (\mathbf{M}_k[l, j])^{-1}\})\}$$

$$\cap \mathrm{UNC}(\{\theta_n\}) \cap \mathrm{LIM}(\{\theta_n\})$$

$$\cap \left\{ \lim_n \mathbf{\Psi}_n \mathbf{M}_n \sum_{k=1}^n (f_k(\theta_0, \nu_0) - f_k(\theta_0, \widehat{\nu})) \mathbf{f}'_k(\theta_0, \widehat{\nu}) \stackrel{P}{=} 0 \right\}$$

$$\cap \left\{ \lim_n \mathbf{\Psi}_n \mathbf{M}_n \sum_{k=1}^n \eta_k \mathbf{f}'_k(\theta_0, \widehat{\nu}) \exists \text{ in distribution} \right\}.$$

*Then, on $\Omega_\infty$,*

(6.1) $$\lim_n \mathbf{\Psi}_n(\widehat{\theta}_n - \theta_0) \stackrel{\mathcal{D}}{=} \lim_n \mathbf{\Psi}_n \mathbf{M}_n \sum_{k=1}^n \eta_k \mathbf{f}'_k(\theta_0, \widehat{\nu}).$$

*Under A2, replace $\widehat{\nu}$ by $\nu$ in $\mathrm{LIP}_\theta$, $\mathrm{VAR}_\theta$ and replace these conditions by $\mathrm{LIP}_{\theta,\nu}$, $\mathrm{VAR}_{\theta,\nu}$, respectively.*

REMARK 6.2. In the linear model $f_n(\theta, \nu) = \theta^T \mathbf{W}_n$ with $\overline{\lim}_n \sigma_n^2 \stackrel{\text{a.s.}}{<} \infty$, we have $\mathbf{M}_n = [\sum_{k=1}^n \mathbf{W}_k \mathbf{W}_k^T]^{-1}$ and, if $\lim_n \widehat{\theta}_n \stackrel{P}{=} \theta_0$, $\Omega_\infty$ is reduced to

$$\Omega_\infty = \left\{ \lim_n \lambda_{\min} \left\{ \sum_{k=1}^n \mathbf{W}_k \mathbf{W}_k^T \right\} \stackrel{\text{a.s.}}{=} \infty \right\}$$

$$\cap \left\{ \lim_n \left[ \sum_{k=1}^n \mathbf{W}_k \mathbf{W}_k^T \right]^{-1/2} \sum_{k=1}^n \eta_k \mathbf{W}_k \exists \text{ in distribution} \right\}.$$

REMARK 6.3. A CLT for martingale arrays may be applied to the right-hand side of (6.1) given $\widehat{\nu}$ only if $\mathbf{f}'_k(\theta_0, \widehat{\nu})$ does not depend on $n$, which is the case under A1.

PROOF OF PROPOSITION 6.1. We derive, as in classical regression, the asymptotic distribution of the estimator from the first order Taylor series expansion of $S'_{n|\widehat{\nu};i}(\widehat{\theta}_n)$ at $\theta_0$ for all $i = 1, \ldots, p$ (see, e.g., [34]):

$$(6.2) \quad S'_{n|\widehat{\nu};i}(\widehat{\theta}_n) = S'_{n|\widehat{\nu};i}(\theta_0) + \sum_{j=1}^{p} [\partial(S'_{n|\widehat{\nu};i})/\partial\theta_j](\theta_n)(\widehat{\theta}_{n,j} - \theta_{0,j}),$$

where $\theta_n = \theta_0 + a_n(\widehat{\theta}_n - \theta_0)$, $a_n \in (0,1)$. Since $S'_{n|\widehat{\nu};i}(\widehat{\theta}_n) = 0$ for all $i$, (6.2) is written in matrix form, $0 = \mathbf{S}'_{n|\widehat{\nu}}(\theta_0) + \mathbf{S}''_{n|\widehat{\nu}}(\theta_n)(\widehat{\theta}_n - \theta_0)$, implying that if $\mathbf{S}''_{n|\widehat{\nu}}(\cdot)$ is invertible in a neighborhood of $\theta_0$, then

$$(6.3) \quad \mathbf{\Psi}_n(\widehat{\theta}_n - \theta_0) = -\mathbf{\Psi}_n[\mathbf{S}''_{n|\widehat{\nu}}(\theta_n)]^{-1}\mathbf{S}'_{n|\widehat{\nu}}(\theta_0)$$

for any $p \times p$ matrix $\mathbf{\Psi}_n$. Moreover, by definition,

$$\mathbf{S}'_{n|\widehat{\nu}}(\theta) = -2\sum_{k=1}^{n}(\eta_k + f_k(\theta_0, \nu_0) - f_k(\theta, \widehat{\nu}))\mathbf{f}'_k(\theta, \widehat{\nu})$$

$$(6.4) \quad \mathbf{S}''_{n|\widehat{\nu}}(\theta) = 2\sum_{k=1}^{n}\mathbf{f}'_k(\theta, \widehat{\nu})\mathbf{f}'^T_k(\theta, \widehat{\nu}) - 2\sum_{k=1}^{n}\eta_k \mathbf{f}''_k(\theta, \widehat{\nu})$$

$$- 2\sum_{k=1}^{n}(f_k(\theta_0, \nu_0) - f_k(\theta, \widehat{\nu}))\mathbf{f}''_k(\theta, \widehat{\nu}).$$

Considering (6.4), if the conditions

$$(6.5) \quad \lim_n \left[\sum_{k=1}^{n}\eta_k \mathbf{f}''_k(\theta_n, \widehat{\nu})\right]\mathbf{M}_n \stackrel{P}{=} 0,$$

$$(6.6) \quad \lim_n \left[\sum_{k=1}^{n}(f_k(\theta_0, \nu_0) - f_k(\theta_n, \widehat{\nu}))\mathbf{f}''_k(\theta_n, \widehat{\nu})\right]\mathbf{M}_n \stackrel{P}{=} 0,$$

$$(6.7) \quad \lim_n \left[\sum_{k=1}^{n}\mathbf{f}'_k(\theta_n, \widehat{\nu})\mathbf{f}'^T_k(\theta_n, \widehat{\nu})\right]\mathbf{M}_n \stackrel{P}{=} \mathbf{I}$$

are fulfilled, then we get $\lim_n 2^{-1}\mathbf{S}''_{n|\widehat{\nu}}(\theta_n)\mathbf{M}_n \stackrel{P}{=} \mathbf{I}$. Therefore, writing

$$\mathbf{\Psi}_n(\widehat{\theta}_n - \theta_0) = -\mathbf{\Psi}_n\mathbf{M}_n[\mathbf{S}''_{n|\widehat{\nu}}(\theta_n)\mathbf{M}_n]^{-1}\mathbf{S}'_{n|\widehat{\nu}}(\theta_0)$$

and using Slutsky's theorem, we get (6.1) on $\Omega_\infty$.

Now, concerning (6.5), it is satisfied if, for all $i,l$,

$$(6.8) \qquad \sup_{l,i,j} \limsup_n \sup_\theta \left| \sum_{k=1}^n \eta_k f''_{k;i,l}(\theta,\widehat{\nu}) \right| |\mathbf{M}_n[l,j]| \overset{\text{a.s.}}{=} 0,$$

which is verified according to Proposition 5.1. Finally, (6.6) and (6.7) are verified on $\text{LIM}(\{\theta_n\}) \cap \text{UNC}(\{\theta_n\})$. $\square$

When $\text{UNC}(\{\theta_n\})$ is not verified, we may instead use the second order Taylor series expansion of $S'_{n|\widehat{\nu};i}(\widehat{\theta}_n)$ at $\theta_0$. So, let $[\partial f''_{k;i,j}/\partial \theta_l](\theta,\widehat{\nu}) =: f'''_{k;i,j,l}(\theta,\widehat{\nu})$ for any $i,j,l$.

PROPOSITION 6.2. *Let us assume that $f_k(\theta,\nu)$ has third derivatives in $\theta$ for each $k$ (and for $\nu = \widehat{\nu}$ under A1 and for $\nu$ under A2) and that there exists a $p \times p$ $\mathcal{F}_{n-1}$-measurable matrix $\mathbf{\Psi}_n$ such that $P(\Omega_\infty) > 0$, where $\Omega_\infty$ is defined under A1 by*

$$\Omega_\infty = \bigcap_{i,j,h,l,m} \left\{ \text{LIP}_\theta(\{f''_{k;i,l}(\theta,\widehat{\nu})\}) \cap \text{SI}(\{(\mathbf{M}_n[l,j])^{-1}\}) \right.$$
$$\cap \text{VAR}_\theta(\{\sigma_k^2, f''_{k,i,l}(\theta,\widehat{\nu}), (\mathbf{M}_k[l,j])^{-1}\}) \cap \text{LIM}(\{\theta_n\})$$
$$\cap \left\{ \overline{\lim}_{n,\theta_n \to \theta_0} \left| (\mathbf{\Psi}_n \mathbf{M}_n)[i,j] \sum_{k=1}^n f'_{k;m}(\theta_n,\widehat{\nu}) f''_{k;h,l}(\theta_n,\widehat{\nu}) \right| \overset{\text{a.s.}}{<} \infty \right\}$$
$$\cap \text{SI}(\{((\mathbf{\Psi} \mathbf{M}_n)[l,j])^{-1}\})$$
$$\cap \text{VAR}_\theta(\{\sigma_k^2, f'''_{k;j,h,l}(\theta,\widehat{\nu}), ((\mathbf{\Psi}_k \mathbf{M}_k)[i,j])^{-1}\})$$
$$\left. \cap \left\{ \overline{\lim}_{n,\theta_n \to \theta_0} (\mathbf{\Psi}_n \mathbf{M}_n)[i,j] \left| \sum_{k=1}^n f_k(\theta_0,\nu_0) f'''_{k;j,h,l}(\theta_n,\widehat{\nu}) \right| \overset{\text{a.s.}}{<} \infty \right\} \right\}$$
$$\cap \left\{ \lim_n \mathbf{\Psi}_n \mathbf{M}_n \sum_{k=1}^n (f_k(\theta_0,\nu_0) - f_k(\theta_0,\widehat{\nu})) \mathbf{f}'_k(\theta_0,\widehat{\nu}) \overset{P}{=} 0 \right\}$$
$$\cap \left\{ \lim_n \mathbf{\Psi}_n \mathbf{M}_n \sum_{k=1}^n \eta_k \mathbf{f}'_k(\theta_0,\widehat{\nu}) \; \exists \; \text{in distribution} \right\}.$$

*Then, on $\Omega_\infty$, $\lim_n \mathbf{\Psi}_n(\widehat{\theta}_n - \theta_0) \overset{\mathcal{D}}{=} \lim_n \mathbf{\Psi}_n \mathbf{M}_n \sum_{k=1}^n \eta_k \mathbf{f}'_k(\theta_0,\widehat{\nu})$.*

*Under A2, replace $\widehat{\nu}$ by $\nu$ in $\text{LIP}_\theta(\cdot)$, $\text{VAR}_\theta(\cdot)$ and replace these conditions by the corresponding conditions on $\theta,\nu$.*

PROOF. Consider the second order Taylor series expansion of $S'_{n|\widehat{\nu};i}(\widehat{\theta}_n)$ at $\theta_0$: for all $i = 1, \ldots, p$,

$$
\begin{aligned}
S'_{n|\widehat{\nu};i}(\widehat{\theta}_n) = {} & S'_{n|\widehat{\nu};i}(\theta_0) + \sum_{j=1}^{p} [\partial(S'_{n|\widehat{\nu};i})/\partial\theta_j](\theta_0)(\widehat{\theta}_{n,j} - \theta_{0,j}) \\
& + \frac{1}{2} \sum_{j,l} [\partial^2(S'_{n|\widehat{\nu};i})/(\partial\theta_j\,\partial\theta_l)](\theta_n)(\widehat{\theta}_{n,l} - \theta_{0,l})(\widehat{\theta}_{n,j} - \theta_{0,j}),
\end{aligned}
\tag{6.9}
$$

where $\theta_n = \theta_0 + a_n(\widehat{\theta}_n - \theta_0)$, $a_n \in \,]0,1[$. Using the definition of $\widehat{\theta}_n$, (6.9) may be written

$$
0 = \mathbf{S}'_{n|\widehat{\nu}}(\theta_0) + \mathbf{S}''_{n|\widehat{\nu}}(\theta_0)(\widehat{\theta}_n - \theta_0) + \frac{1}{2}\sum_{l=1}^{p}[\partial(\mathbf{S}''_{n|\widehat{\nu}})/\partial\theta_l](\theta_n)(\widehat{\theta}_{n,l} - \theta_{0,l})(\widehat{\theta}_n - \theta_0).
$$

Let $\boldsymbol{\Psi}_n$ a $p \times p$ matrix. Then

$$
\boldsymbol{\Psi}_n\!\left(\mathbf{I} + [\mathbf{S}''_{n|\widehat{\nu}}(\theta_0)]^{-1}\!\left[\frac{1}{2}\sum_l [\partial(\mathbf{S}''_{n|\widehat{\nu}})/\partial\theta_l](\theta_n)(\widehat{\theta}_{n,l} - \theta_{0,l})\right]\right)(\widehat{\theta}_n - \theta_0)
$$
$$
= -\boldsymbol{\Psi}_n[\mathbf{S}''_{n|\widehat{\nu}}(\theta_0)]^{-1}\mathbf{S}'_{n|\widehat{\nu}}(\theta_0).
$$

The proof is then similar to the proof of Proposition 6.1. To prove that $\lim_n \mathbf{S}''_{n|\widehat{\nu}}(\theta_0)\mathbf{M}_n \overset{P}{=} 2I$, we need the assumptions of Proposition 6.1, where $\{\theta_n\}$ and $\theta$ are replaced by $\theta_0$. Next, using the fact that $\lim_n \mathbf{S}''_{n|\widehat{\nu}}(\theta_0)\mathbf{M}_n \overset{P}{=} 2I$, we show that

$$
\lim_n \boldsymbol{\Psi}_n\!\left([\mathbf{S}''_{n|\widehat{\nu}}(\theta_0)]^{-1}\!\left[\frac{1}{2}\sum_l [\partial(\mathbf{S}''_{n|\widehat{\nu}})/\partial\theta_l](\theta_n)(\widehat{\theta}_{n,l} - \theta_{0,l})\right]\right)(\widehat{\theta}_n - \theta_0) \overset{P}{=} 0.
$$

From (6.4), we deduce that

$$
\begin{aligned}
\frac{1}{2}\frac{\partial \mathbf{S}''_{n|\widehat{\nu}}}{\partial\theta_l}(\theta)[j,h] = {} & \sum_{k=1}^{n} f''_{k;j,l}(\theta,\widehat{\nu})f'_{k;h}(\theta,\widehat{\nu}) + \sum_{k=1}^{n} f''_{k;h,l}(\theta,\widehat{\nu})f'_{k;j}(\theta,\widehat{\nu}) \\
& + \sum_{k=1}^{n} f''_{k;j,h}(\theta,\widehat{\nu})f'_{k;l}(\theta,\widehat{\nu}) - \sum_{k=1}^{n} f_k(\theta_0,\nu_0)f'''_{k;j,h,l}(\theta) \\
& - \sum_{k=1}^{n} \eta_k f'''_{k;j,h,l}(\theta,\widehat{\nu}).
\end{aligned}
$$

The proof is then as in Proposition 6.1. $\square$

**7. Estimation of $\nu_0$.** In the previous sections, we obtained conditions leading to the consistency and the asymptotic distribution of $\widehat{\theta}_n$ given $\{\widehat{\nu}_n\}$. When the asymptotic distribution of $\widehat{\theta}_n$ depends on $\nu_0$, it is necessary to derive a consistent estimator of $\nu_0$. Let us write $\nu_0 = (\nu_0^{(1)}, \nu_0^{(2)})$, where $\nu_0^{(1)}$ is of dimension $q_1 \leq \infty$ and $\nu_0^{(2)}$ is of dimension $q_2 < \infty$, and such that the asymptotic distribution of $\widehat{\theta}_n$ is independent of $\nu_0^{(1)}$ and

$$f_n(\theta_0, \nu_0) = f_n(\theta_0, \nu_0^{(1)}), \qquad \sigma_n^2(\theta_0, \nu_0) = \sigma_n^2(\theta_0, \nu_0^{(2)}),$$

where, now, $\sigma_n^2(\theta_0, \nu_0^{(2)}) := E([Y_n - E(Y_n|\mathcal{F}_{n-1})]^2|\mathcal{F}_{n-1})$. Let us further assume that there exists $s(\theta, \nu^{(2)})$, a continuous function in $(\theta, \nu^{(2)})$, such that $\lim_n \Psi_n(\widehat{\theta}_n - \theta_0)[s(\theta_0, \nu_0^{(2)})]^{-1} \stackrel{\mathcal{D}}{=} \mathcal{L}$, where $\mathcal{L}$ is independent of the unknown parameters. Then, for any sequence of estimators $\{(\widehat{\theta}_n, \widehat{\nu}_n^{(2)})\}$ such that $\lim_n(\widehat{\theta}_n, \widehat{\nu}_n^{(2)}) \stackrel{P}{=} (\theta_0, \nu_0^{(2)})$, thanks to Slutsky's theorem, $\lim_n \Psi_n(\widehat{\theta}_n - \theta_0)[s(\widehat{\theta}_n, \widehat{\nu}_n^{(2)})]^{-1} \stackrel{\mathcal{D}}{=} \mathcal{L}$, which allows for the elaboration of confidence regions of $\theta_0$.

Therefore, we build here a CLSE of $\nu_0^{(2)}$ and give the conditions for its consistency. Let us define

$$\widehat{\nu}_n^{(2)} = \arg\min_{\nu^{(2)} \in \mathcal{N}^{(2)}} \widetilde{S}_{n|\widehat{\theta}_n, \widehat{\nu}^{(1)}}(\nu^{(2)}),$$

$$\widetilde{S}_{n|\widehat{\theta}_n, \widehat{\nu}^{(1)}}(\nu^{(2)}) := \sum_{k=1}^{n}(\widetilde{Y}_{k,n} - \widetilde{f}_{k,n}(\widehat{\theta}_n, \nu^{(2)}))^2,$$

where $\widetilde{Y}_{k,n} = (Y_k - f_k(\widehat{\theta}_n, \widehat{\nu}^{(1)}))^2 = (\eta_k + d_k(\widehat{\theta}_n, \widehat{\nu}^{(1)}))^2$, $\widetilde{f}_{k,n}(\theta_0, \nu_0^{(2)}) := E(\widetilde{Y}_{k,n}|\mathcal{F}_{k-1}) = \sigma_k^2(\theta_0, \nu_0^{(2)}) + d_k^2(\widehat{\theta}_n, \widehat{\nu}^{(1)})$. Let us define

$$\widetilde{d}_k(\nu^{(2)}|\widehat{\theta}_n) := \widetilde{f}_{k,n}(\widehat{\theta}_n, \nu_0^{(2)}) - \widetilde{f}_{k,n}(\widehat{\theta}_n, \nu^{(2)}) := \sigma_k^2(\widehat{\theta}_n, \nu_0^{(2)}) - \sigma_k^2(\widehat{\theta}_n, \nu^{(2)}),$$

$$\widetilde{d}_k(\widehat{\theta}_n|\nu^{(2)}) := \widetilde{f}_{k,n}(\theta_0, \nu^{(2)}) - \widetilde{f}_{k,n}(\widehat{\theta}_n, \nu^{(2)}) := \sigma_k^2(\theta_0, \nu^{(2)}) - \sigma_k^2(\widehat{\theta}_n, \nu^{(2)}),$$

$$\widetilde{D}_n(\nu^{(2)}|\widehat{\theta}_n) := \sum_{k=1}^{n} \widetilde{d}_k^2(\nu^{(2)}|\widehat{\theta}_n), \qquad \widetilde{D}_n(\widehat{\theta}_n|\nu^{(2)}) := \sum_{k=1}^{n} \widetilde{d}_k^2(\widehat{\theta}_n|\nu^{(2)}),$$

$$\widetilde{L}_{n,n}(\nu^{(2)}|\widehat{\theta}_n) := \sum_{k=1}^{n} \widetilde{\eta}_{k,n} \widetilde{d}_k(\nu^{(2)}|\widehat{\theta}_n),$$

$$\widetilde{\eta}_{k,n} := \widetilde{Y}_{k,n} - E(\widetilde{Y}_{k,n}|\mathcal{F}_{k-1}) = \widetilde{\eta}_k + 2\eta_k d_k(\widehat{\theta}_n, \widehat{\nu}^{(1)}),$$

where $\widetilde{\eta}_k := \eta_k^2 - \sigma_k^2(\theta_0, \nu_0^{(2)})$. This implies that

(7.1) $\qquad \widetilde{L}_{n,n}(\nu^{(2)}|\widehat{\theta}_n) = \widetilde{L}_n(\nu^{(2)}|\widehat{\theta}_n) + 2L_n(\nu^{(2)}|\widehat{\theta}_n, \widehat{\nu}^{(1)}),$

where
$$\widetilde{L}_n(\nu^{(2)}|\widehat{\theta}_n) := \sum_{k=1}^{n} \widetilde{\eta}_k \widetilde{d}_k(\nu^{(2)}|\widehat{\theta}_n),$$

$$L_n(\nu^{(2)}|\widehat{\theta}_n, \widehat{\nu}^{(1)}) := \sum_{k=1}^{n} \eta_k d_k(\widehat{\theta}_n, \widehat{\nu}^{(1)}) \widetilde{d}_k(\nu^{(2)}|\widehat{\theta}_n).$$

Let $\widetilde{\sigma}_k^2 := E(\widetilde{\eta}_k^2 | \mathcal{F}_{k-1})$.

PROPOSITION 7.1. *Let us assume that* $\sigma_k^2(\theta, \nu_0^{(2)})$ *is continuous in* $\theta$, *that* $E(\eta_k^4|\mathcal{F}_{k-1}) \stackrel{a.s.}{<} \infty$ *for all* $k$ *and that* $\lim_n \widehat{\theta}_n \stackrel{a.s.}{=} \theta_0$. *Let* $\Omega_\infty \subset \Omega$ *be defined by*

$$\Omega_\infty = \mathrm{LIP}_{\theta,\nu^{(2)}}(\{\widetilde{d}_k(\nu^{(2)}|\theta)\}) \cap \mathrm{LIP}_{\theta,\nu^{(2)}}(\{d_k(\theta,\widehat{\nu}^{(1)})\widetilde{d}_k(\nu^{(2)}|\theta)\})$$

$$\cap \mathrm{SI}_{\theta,\nu^{(2)}}(\{\widetilde{D}_n(\nu^{(2)}|\theta)\}) \cap \mathrm{VAR}_{\theta,\nu^{(2)}}(\{\widetilde{\sigma}_k^2, \widetilde{d}_k(\nu^{(2)}|\theta), \widetilde{D}_k(\nu^{(2)}|\theta)\})$$

$$\cap \mathrm{VAR}_{\theta,\nu^{(2)}}(\{\sigma_k^2, d_k(\theta,\widehat{\nu}^{(1)})\widetilde{d}_k(\nu^{(2)}|\theta), \widetilde{D}_k(\nu^{(2)}|\theta)\}).$$

*Let us assume that* $P(\Omega_\infty) > 0$. *Then* $\lim_n \widehat{\nu}_n^{(2)} \stackrel{a.s.}{=} \nu_0^{(2)}$ *on* $\Omega_\infty$.

PROOF. Let us write
$$\widetilde{Y}_{k,n} - \widetilde{f}_k(\widehat{\theta}_n, \nu^{(2)}) = \widetilde{\eta}_{k,n} + \widetilde{d}_k(\widehat{\theta}_n|\nu_0^{(2)}) + \widetilde{d}_k(\nu^{(2)}|\widehat{\theta}_n).$$

Then, using Hölder's inequality, this decomposition leads to
$$\widetilde{S}_{n|\widehat{\theta}_n, \widehat{\nu}^{(1)}}(\nu^{(2)}) - \widetilde{S}_{n|\widehat{\theta}_n, \widehat{\nu}^{(1)}}(\nu_0^{(2)})$$

$$= \sum_{k=1}^{n}[\widetilde{d}_k(\nu^{(2)}|\widehat{\theta}_n)]^2 + 2\sum_{k=1}^{n}(\widetilde{\eta}_{k,n} + \widetilde{d}_k(\widehat{\theta}_n|\nu_0^{(2)}))\widetilde{d}_k(\nu^{(2)}|\widehat{\theta}_n)$$

$$\geq \widetilde{D}_n(\nu^{(2)}|\widehat{\theta}_n) - 2|\widetilde{L}_{n,n}(\nu^{(2)}|\widehat{\theta}_n)| - 2[\widetilde{D}_n(\widehat{\theta}_n|\nu_0^{(2)})]^{1/2}[\widetilde{D}_n(\nu^{(2)}|\widehat{\theta}_n)]^{1/2},$$

which implies that
$$\inf_{\nu^{(2)} \in B_\delta^c} \widetilde{S}_{n|\widehat{\theta}_n, \widehat{\nu}^{(1)}}(\nu^{(2)}) - \widetilde{S}_{n|\widehat{\theta}_n, \widehat{\nu}^{(1)}}(\nu_0^{(2)})$$

$$\geq \inf_{\nu^{(2)} \in B_\delta^c} \widetilde{D}_n(\nu^{(2)}|\widehat{\theta}_n)$$

$$\times \left[1 - 2 \sup_{(\theta,\nu^{(2)}) \in B_\delta^c} \frac{|\widetilde{L}_{n,n}(\nu^{(2)}|\theta)|}{\widetilde{D}_n(\nu^{(2)}|\theta)} - 2\frac{[\widetilde{D}_n(\widehat{\theta}_n|\nu_0^{(2)})]^{1/2}}{\inf_{\nu^{(2)} \in B_\delta^c}[\widetilde{D}_n(\nu^{(2)}|\widehat{\theta}_n)]^{1/2}}\right].$$

The result then follows from the assumptions of the proposition, from (7.1), from the SLLNSM applied to $\sup_{(\theta,\nu^{(2)}) \in B_\delta^c} |\widetilde{L}_n(\nu^{(2)}|\theta)|[\widetilde{D}_n(\nu^{(2)}|\theta)]^{-1}$ and to $\sup_{(\theta,\nu^{(2)})\in B^c_\delta}|L_n(\nu^{(2)}|\theta,\widehat{\nu}^{(1)})|[\widetilde{D}_n(\nu^{(2)}|\theta)]^{-1}$, and from $\lim_n|\widetilde{d}_n^2(\widehat{\theta}_n|\nu_0^{(2)})| \times [\inf_{\nu^{(2)}\in B^c_\delta}|\widetilde{d}_n^2(\nu^{(2)}|\widehat{\theta}_n)|]^{-1} 1_{\{\inf_{\nu^{(2)}\in B^c_\delta}\widetilde{d}_n(\nu^{(2)}|\widehat{\theta}_n)\neq 0\}} \stackrel{a.s.}{=} 0$, which is deduced from the continuity of $\sigma_k^2(\theta,\nu_0^{(2)})$ in $\theta$ and which, according to Toeplitz's lemma, leads to $\lim_n \widetilde{D}_n(\widehat{\theta}_n|\nu_0^{(2)})[\inf_{\nu^{(2)}\in B^c_\delta}\widetilde{D}_n(\nu^{(2)}|\widehat{\theta}_n)]^{-1} \stackrel{a.s.}{=} 0$. □

**8. Examples in branching processes.** Here, we deal with single-type Markovian branching processes in discrete time. The process $\{N_n\}$ of population sizes at each time is defined by $N_n = \sum_{i=1}^{N_{n-1}} X_{n,i}$, where the $\{X_{n,i}\}_i$ are i.i.d. given $\mathcal{F}_{n-1}$ that is generated by $\{N_k\}_{k\leq n-1}$. In each example, the normalization is carried out in order to get $\overline{\lim}_n E(\eta_n^2|\mathcal{F}_{n-1}) \stackrel{a.s.}{<} \infty$, implying that assumption $\text{VAR}_\theta(\cdot)$ of Proposition 5.1 is satisfied.

8.1. *BGW process.* The $\{X_{n,i}\}_i$ are i.i.d. $(m_0,\sigma_0^2)$ and are independent of $\mathcal{F}_{n-1}$. Then $M_Y$ is defined by

$$Y_n = N_n N_{n-1}^{-1/2}, \qquad f_n(m_0) = m_0 N_{n-1}^{1/2}, \qquad \eta_n = N_{n-1}^{-1/2} \sum_{i=1}^{N_{n-1}} (X_{n,i} - m_0).$$

Therefore, $E(\eta_n^2|\mathcal{F}_{n-1}) = \sigma_0^2$. Let $\theta_0 = m_0$. Then $f_n(m)$ is Lipschitz in $m$ with $g_n = N_{n-1}^{1/2}$ and $\inf_{m\in B^c_\delta} D_n(m) = \inf_{m\in B^c_\delta}(m_0 - m)^2 \sum_{k=1}^n N_{k-1}$, which converges a.s. to $\infty$, as $n \to \infty$, on the nonextinction set $\Omega_\infty$, where $P(\Omega_\infty) > 0$, for $m_0 > 1$. Therefore, assuming that $m_0 > 1$ and using Proposition 3.1, we have $\lim_n \widehat{m}_n \stackrel{a.s.}{=} m_0$ on $\Omega_\infty$.

REMARK 8.1 (Direct proof). Recall that $\widehat{m}_n$ is also the MLE estimator of $m_0$ [6] and the Harris estimator [11]:

$$\widehat{m}_n = \frac{\sum_{k=1}^n N_k}{\sum_{k=1}^n N_{k-1}} = \frac{\sum_{k=1}^n (N_k m_0^{-k}) m_0^k}{\sum_{k=1}^n m_0^k} \frac{m_0 \sum_{k=1}^n m_0^{k-1}}{\sum_{k=1}^n (N_{k-1} m_0^{-(k-1)}) m_0^{k-1}}.$$

The strong consistency of $\{\widehat{m}_n\}$ on $\Omega_\infty$ is then obtained classically (see [8]) by using Toeplitz's lemma with $\lim_n N_n m_0^{-n} \stackrel{a.s.}{=} W_{N_0}$, where $W_{N_0} \stackrel{a.s.}{>} 0$ on $\Omega_\infty$, $E(W_{N_0}) = N_0$. Note that the indirect proof based on Proposition 3.1 does not require knowledge of the asymptotic behavior of the process as is the case in this direct proof.

In this model, (6.1) becomes

$$(8.1) \quad \lim_n \Psi_n(\widehat{m}_n - m_0) = \lim_n \Psi_n \left[\sum_{k=1}^n N_{k-1}\right]^{-1} \sum_{k=1}^n \sum_{i=1}^{N_{k-1}} (X_{k,i} - m_0),$$

which is also the expression obtained classically by directly using the expression for $\widehat{m}_n$. Since $\sum_{k=1}^{n}[f'_k(m_0)]^2 = \sum_{k=1}^{n} N_{k-1} =: S_{n-1}$, $\sum_{k=1}^{n} \eta_k f'_k(m_0) = \sum_{k=1}^{n} \sum_{i=1}^{N_{k-1}} (X_{k,i} - m_0) =: \sum_{j=1}^{S_{n-1}} (X_j - m_0)$, where the $\{X_j\}$ correspond a.s. to the $\{X_{k,i}\}_{i,k}$, ordered according to $i$ and then $k$. Setting $\Psi_n = [\sum_{k=1}^{n} m_0^{k-1}]^{1/2}$, (8.1) then becomes

$$\Psi_n(\widehat{m}_n - m_0) = \left(\left[\sum_{k=1}^{n} m_0^{k-1}\right]^{1/2} \left[\sum_{k=1}^{n} N_{k-1}\right]^{-1/2}\right) \left(S_{n-1}^{-1/2} \sum_{j=1}^{S_{n-1}} (X_j - m_0)\right),$$

which leads, thanks to Toeplitz's lemma on the first term of the right-hand side and a CLT for random sums [2] on the second term, to $\lim_n \Psi_n(\widehat{m}_n - m_0) \stackrel{\mathcal{D}}{=} W_{N_0}^{-1} U$, where $U \sim \mathcal{N}(0, \sigma_0^2)$, $U$ and $W_{N_0}$ being independent.

Similar results may be obtained for $\widehat{m}_{h,n}$, where $n - h$ is constant. In this case, "$\sum_{k=1}^{n}$" must be replaced by "$\sum_{k=h+1}^{n}$." When $h = n-1$, $\widehat{m}_{h,n} = N_n/N_{n-1}$ is the Lotka–Nagaev estimator.

Since $m_0 = 1$ is a threshold for the asymptotic behavior of $\{N_n\}$, if we do not know a priori whether $m < 1$ or $m > 1$, then we may estimate $m$ using $E(N_n|N_{n-1}, N_n \neq 0)$ instead of $E(N_n|N_{n-1})$.

8.2. *Size-dependent branching processes.* $N_n = \sum_{i=1}^{N_{n-1}} X_{n,i}$ with $\{X_{n,i}\}$ i.i.d. $(m_{\theta_0}(N_{n-1}), \sigma_0^2(N_{n-1}))$, $\lim_N m_{\theta_0}(N) = m_0$, $\sigma_0^2(N) = O(N^{\beta_0})$, where $\beta_0$ is assumed to be known [19]. The model $M_Y$ is then

$$Y_n = N_n N_{n-1}^{-(1+\beta_0)/2} = m_{\theta_0}(N_{n-1}) N_{n-1}^{(1-\beta_0)/2} + \eta_n.$$

A particular example of a size-dependent branching process is the process modeling the amplification process in the polymerase chain reaction setting, taking into account the saturation phenomenon due to the closed medium [17, 30]. The ultimate goal of this technology is the estimation of $N_0$, the initial number of DNA molecule fragments, through the amplified population in vitro. At each cycle of the amplification process, a DNA fragment may product two DNA fragments by replication after three successive steps—heating, annealing and synthesis. The amplification process exhibits three different phases: the exponential phase, during which the replication is not limited, then a saturation phase involving a "linear" phase, followed by a "plateau" phase where the replication is less and less efficient. Since the size of the population increases very quickly (exponential increase during the first cycles) and since the observation errors are very important during the first cycles but become more and more negligible relative to the signal as the number of cycles increases, the estimations should be based on the observations starting only from the end of the exponential phase. The classical estimation method is based on a regression model using the observations at the end of the exponential phase of a set of amplification processes starting from successive dilutions of a given DNA sample and assumed to have the same replication probability. Besides the probable violation of the basic assumptions of this method (i.i.d. errors, identical replication probabilities, identical initial size up to the dilution factor, etc.), this method is costly since the accuracy of the estimator requires a large number of such trajectories. Therefore, we developed conditional least squares estimation based on a single amplification process (or two if the "observation unit to number of molecules" conversion is required) [13, 24, 27]. Moreover, because of a very large population after only a few cycles, the asymptotic properties are obtained from the end of the exponential phase, leading to a great accuracy of the estimator of the replication probability, which is crucial for the accuracy of the estimation of $N_0$. We focus here on the estimation of the replication probability based on different models.

The amplification process may be modeled by a simple branching process $N_n = \sum_{i=1}^{N_{n-1}} X_{n,i}$, where the $\{X_{n,i}\}$ are i.i.d. given $\mathcal{F}_{n-1}$ with $P(X_{n,i} = 1|\mathcal{F}_{n-1}) = 1 - p_n$, $p_n := P(X_{n,i} = 2|\mathcal{F}_{n-1})$ being the probability of replication at the $n$th cycle. When restricting the modeling to the exponential phase, we may assume that the $\{X_{n,i}\}$ are i.i.d. $(m_0, \sigma_0^2)$, where $m_0 = 1 + p_0$, $p_0 = p_n$ and $\sigma_0^2 = p_0(1 - p_0)$, that is, the process is a BGW branching process ([13, 27], see previous item).

We now take into account the saturation phase. Therefore, the probability of replication is a decreasing function of the current size of the population. For example, Schnell and Mendoza [30] proposed the following enzymological model:

$$\text{M1: } p_n = [K_0]([K_0] + [N_{n-1}])^{-1} = K_0(K_0 + N_{n-1})^{-1},$$

where $[K_0]$ is the Michaelis–Menten constant, $[N_{n-1}]$ is the concentration of molecules at time $n - 1$ and $K_0 = [K_0] \times V$, where $V$ is the volume of the reaction. Then $\{N_n\}$ is a near-critical process with an a.s. nonextinction and $\lim_n N_n n^{-1} \stackrel{\text{a.s.}}{=} K_0$ [17].

This model may be generalized in order to take into account a saturation threshold $S_0 \geq N_0$. For example, in [24],

$$\text{M2: } p_n = \frac{K_0}{K_0 + N_{S_0,n-1}} \frac{[1 + \exp(-C_0(S_0^{-1} N_{S_0,n-1} - 1))]}{2},$$

where $N_{S_0,n-1} = S_0$ if $N_{n-1} < S_0$, and $N_{S_0,n-1} = N_{n-1}$ if $N_{n-1} \geq S_0$. When $C_0 = 0$, $S_0 = N_0$, M2 is reduced to M1. Since M2 is a BGW process for all $n$ such that $N_{n-1} < S_0$, it follows that when $S_0 \to \infty$, $M_2$ tends to a BGW process. In the general case $C_0 \neq 0$ with $S_0 < \infty$, we have $\lim_n N_n n^{-1} = K_0/2$. When setting $\theta = (K, C, S)$, we have $\lim_n \inf_\theta D_n(\theta) \stackrel{\text{a.s.}}{<} \infty$ and $K_0$ is the only parameter verifying $\lim_n \inf_K D_n(K) \stackrel{\text{a.s.}}{=} \infty$. The asymptotic properties of the CLSE of $K_0$, assuming $(C_0, S_0) =: \nu_0$, may be found in [24]. They can also be derived in the same way as for the following model using the general results developed here.

Let us now assume another generalization of the Schnell–Mendoza model:

$$\text{M3: } p_n = \left(\frac{K_0}{K_0 + N_{S_0,n-1}}\right)\frac{(1 + S_0^{\alpha_0} N_{S_0,n-1}^{-\alpha_0})}{2}, \qquad \alpha_0 > 0.$$

As for M2, the limit of this model, as $S_0 \to \infty$, is reduced to the BGW model and when $\alpha_0 = 0$ with $S_0 = N_0$, the model is reduced to M1. Let us assume here that $S_0 < \infty$ and $0 < \alpha_0 < \infty$. As in [17] and [24], the asymptotic behavior of the process is linear: $\lim_n N_n n^{-1} \stackrel{a.s.}{=} K_0/2$ and the stochastic regression model is $Y_n = N_n = f_n(\theta_0, \nu_0) + \eta_n$, where

$$f_n(\theta_0, \nu_0) = \left[1 + \left(\frac{K_0}{K_0 + N_{S_0,n-1}}\right)\frac{(1 + S_0^{\alpha_0} N_{S_0,n-1}^{-\alpha_0})}{2}\right] N_{n-1},$$

(8.2)
$$\eta_n = \sum_{i=1}^{N_{n-1}} (X_{n,i} - E(X_{n,i}|N_{n-1})),$$

$$\sigma_n^2 = N_{n-1} p_n (1 - p_n) = O(1).$$

Let $\theta = (K, S^\alpha, \alpha)$, $q = 0$. Using the first order Taylor series development, we have

$$\inf_{\theta \in B_c^\delta} D_n(\theta) = \inf_{|\alpha - \alpha_0| \geq \delta} O\left(\sum_{k=1}^n (k^{-\alpha_0} - k^{-\alpha})^2\right)$$

$$= O\left(\inf_{\widetilde\alpha = \alpha_0 + a_{\delta,\alpha}(\alpha - \alpha_0)} \sum_{k=1}^n [\ln(k)]^2 k^{-2\widetilde\alpha}\right), \qquad a_{\delta,\alpha} \in \,]0,1[,$$

which converges a.s. to $\infty$ for $0 \leq 2\widetilde\alpha \leq 1$. Therefore, if we assume that $\alpha_0 \in \,]0, 1/2[$, then $\lim_n (\widehat K_n, \widehat{S_n^\alpha}, \widehat\alpha_n) \stackrel{a.s.}{=} (K_0, S_0^{\alpha_0}, \alpha_0)$. Next, since the second derivative of $f_n(\theta)$ in $\alpha$ increases to $\infty$ with $n$ more quickly than its first derivative, assumptions such that $\text{VAR}_\theta(\{\sigma_k^2, f''_{k;1,2}(\theta, \widehat\nu), (\mathbf{M}_k[2,1])^{-1}\})$ cannot be verified. Therefore, we restrict the study of the asymptotic distribution to $\widehat K_n$ and then to $(\widehat K_n, \widehat{S_n^\alpha})$.

So, let us first assume that $\theta_0 = K_0$ and $(S_0^{\alpha_0}, \alpha_0) =: \nu_0$ are nuisance parameters estimated by $(\widehat{S_{n_0}^\alpha}, \widehat\alpha_{n_0})$. We have, for $N_{k-1}$ large enough,

$$f'_k(K_0, \widehat\nu) = \frac{N_{k-1}^2}{(K_0 + N_{k-1})^2}\left(\frac{1 + \widehat{S_{n_0}^\alpha} N_{k-1}^{-\widehat\alpha_{n_0}}}{2}\right) = O(1),$$

$$f''_k(K_0, \widehat\nu) = -\frac{N_{k-1}^2}{(K_0 + N_{k-1})^3}(1 + \widehat{S_{n_0}^\alpha} N_{k-1}^{-\widehat\alpha_{n_0}}) = O(k^{-1}).$$

Let us define

$$\Phi_n^2 = \sum_{k=1}^{n}[(1+\widehat{S_{n_0}^{\alpha}}2^{\widehat{\alpha}_{n_0}}K_0^{-\widehat{\alpha}_{n_0}}(k-1)^{-\widehat{\alpha}_{n_0}})/2]^2, \qquad \Psi_n = \Phi_n^{-1}\sum_{k=1}^{n}[f_k'(K_0,\widehat{\nu})]^2.$$

Then, thanks to Toeplitz's lemma, $\lim_n \Phi_n^{-1}\Psi_n \stackrel{\text{a.s.}}{=} 1$ and $\Psi_n = O(1)$. Consequently, for $0 < 2\widehat{\alpha}_{n_0} < 1$, the assumptions of Proposition 6.1, except the last one, which we consider now, are verified. By the definition of $\Psi_n$,

$$(8.3) \quad \Psi_n\left[\sum_{k=1}^{n}[f_k'(K_0,\widehat{\nu})]^2\right]^{-1}\sum_{k=1}^{n}\eta_k f_k'(K_0,\widehat{\nu}) = \Phi_n^{-1}\sum_{k=1}^{n}\eta_k f_k'(K_0,\widehat{\nu}).$$

Moreover,

$$\lim_{k\leq n, k\to\infty} \mathbf{f}_k'(K_0,\widehat{\nu}) \stackrel{\text{a.s.}}{=} 1/2, \qquad \lim_k E(\eta_k^2|\mathcal{F}_{k-1}) \stackrel{\text{a.s.}}{=} K_0/2.$$

Therefore, we may derive the asymptotic distribution of the right-hand side of (8.3) by a classical CLT for martingale arrays ([4, 29]: "let $\{\mathbf{M}_k^{(n)}\}_{k\leq n}$ be a multidimensional martingale triangular array $\sum_{l=1}^{k}\xi_{l,n}$, $E(\xi_{l,n}|\mathcal{F}_{l-1}^{(n)}) \stackrel{\text{a.s.}}{=} 0$. Let us assume (a) $\lim_n \langle\mathbf{M}\rangle_n^{(n)} \stackrel{P}{=} \Gamma$, where $\Gamma$ is a semi-definite deterministic matrix and $\langle\mathbf{M}\rangle_n^{(n)} := \sum_{l=1}^{n}E(\xi_{l,n}\xi_{l,n}^T|\mathcal{F}_{l-1}^{(n)})$; (b) for all $\epsilon > 0$, $\lim_n \sum_{k=1}^{n}E(\|\xi_{k,n}\|^2 1_{\|\xi_{k,n}\|\geq\epsilon}|\mathcal{F}_{k-1}^{(n)}) \stackrel{P}{=} 0$. Then $\lim_n \mathbf{M}_n^{(n)} \stackrel{\mathcal{D}}{=} \mathcal{N}(0,\Gamma)$.") So, let us define the martingale array $\mathbf{M}_k^{(n)} = \Phi_n^{-1}\sum_{l=1}^{k}\eta_l f_l'(K_0,\widehat{\nu})$. Then

$$\langle\mathbf{M}\rangle_n^{(n)} = \Phi_n^{-2}\sum_{k=1}^{n}E(\eta_k^2|\mathcal{F}_{k-1})[f_k'(K_0,\widehat{\nu})]^2,$$

which implies that $\lim_n \langle\mathbf{M}\rangle_n^{(n)} \stackrel{\text{a.s.}}{=} \lim_n E(\eta_n^2|\mathcal{F}_{n-1}) = K_0/2$. Moreover, for $n$ large, using Hölder's and Markov's inequalities, we have

$$\sum_{k=1}^{n}E(\|\xi_{k,n}\|^2 1_{\|\xi_{k,n}\|\geq\epsilon}|\mathcal{F}_{k-1}) \leq \Phi_n^{-2}\sum_{k=1}^{n}E(\eta_k^2 1_{\{\eta_k^2\geq\epsilon^2\Phi_n^2\}}|\mathcal{F}_{k-1})$$

$$\leq [\epsilon\Phi_n^3]^{-1}\sum_{k=1}^{n}[E(\eta_k^4|\mathcal{F}_{k-1})]^{1/2}\sigma_k.$$

Using (8.2), we get $E(\eta_k^4|\mathcal{F}_{k-1}) = O(1)$ and $\sigma_k = O(1)$, which implies the Lindeberg condition since $\Phi_n^2 = O(n)$. Consequently, $\lim_n \Phi_n(\widehat{K}_n - K_0) \stackrel{\mathcal{D}}{=} \mathcal{N}(0, K_0/2)$ or, equivalently, $\lim_n \sqrt{n}(\widehat{K}_n - K_0)(2K_0)^{-1/2} \stackrel{\mathcal{D}}{=} \mathcal{N}(0,1)$. Note that if $\alpha_0 = 0$ (the Schnell–Mendoza model), then, in the same way, $\lim_n \sqrt{n}(\widehat{K}_n - K_0)K_0^{-1/2} \stackrel{\mathcal{D}}{=} \mathcal{N}(0,1)$.

We may also derive the asymptotic distribution using the CLT for random sums [28] applied to the right-hand side of (8.3) using (8.2). According to this CLT, the asymptotic behavior is obtained, even for small $n$, because of the expression of $\eta_k$ as the sum of a large number of centered variables. The asymptotic distribution of the estimator allows for the derivation of accurate confidence intervals of $K_0$.

Let us now assume that $\theta = (K, S^{\alpha_0})$ with $\alpha_0 = \nu_0$ estimated by $\widehat{\alpha}_{n_0}$. We have, for $k \leq n$,

$$f'_{k;1}(\theta,\widehat{\nu}) = N_{k-1}^2(K + N_{k-1})^{-2}(1 + S^{\alpha_0}N_{k-1}^{-\widehat{\alpha}_{n_0}})2^{-1} = O(1),$$

$$f'_{k;2}(\theta,\widehat{\nu}) = KN_{k-1}(K + N_{k-1})^{-1}N_{k-1}^{-\widehat{\alpha}_{n_0}}2^{-1} = O(k^{-\widehat{\alpha}_{n_0}}),$$

$$f''_{k;1,1}(\theta,\widehat{\nu}) = -N_{k-1}^2(K + N_{k-1})^{-3}(1 + S^{\alpha_0}N_{k-1}^{-\widehat{\alpha}_{n_0}}) = O(k^{-1}),$$

$$f''_{k;2,2}(\theta,\widehat{\nu}) = 0,$$

$$f''_{k;1,2}(\theta,\widehat{\nu}) = N_{k-1}^2(K + N_{k-1})^{-2}N_{k-1}^{-\widehat{\alpha}_{n_0}}2^{-1} = f''_{k;2,1}(\theta,\widehat{\nu}) = O(k^{-\widehat{\alpha}_{n_0}}).$$

Therefore, UNC($\{\theta_n\}$) is checked and

$$4\sum_{k=1}^{n} \mathbf{f}'_k(\theta,\widehat{\nu})\mathbf{f}'^T_k(\theta,\widehat{\nu}) = O\begin{pmatrix} n & Kn^{1-\widehat{\alpha}_{n_0}} \\ Kn^{1-\widehat{\alpha}_{n_0}} & K^2 a_n(\widehat{\alpha}_{n_0}) \end{pmatrix},$$

where $a_n(\widehat{\alpha}_{n_0}) = n^{1-2\widehat{\alpha}_{n_0}}1_{\{2\widehat{\alpha}_{n_0}<1\}} + \ln n 1_{\{2\widehat{\alpha}_{n_0}=1\}}$, implying that

$$\mathbf{M}_n = \begin{pmatrix} O(n^{-1}) & O(n^{-\widehat{\alpha}_{n_0}}[a_n(\widehat{\alpha}_{n_0})]^{-1}) \\ O(n^{-\widehat{\alpha}_{n_0}}[a_n(\widehat{\alpha}_{n_0})]^{-1}) & O([a_n(\widehat{\alpha}_{n_0})]^{-1}) \end{pmatrix}.$$

Consequently, all of the assumptions of Proposition 6.1 are satisfied for $\{0 < 2\widehat{\alpha}_{n_0} < 1\}$. Let us define $\mathbf{\Psi}_n = \Phi_n^{-1}[\sum_{k=1}^{n} \mathbf{f}'_k(\theta_0,\widehat{\nu})\mathbf{f}'^T_k(\theta_0,\widehat{\nu})]$, where

$$\Phi_n^2 = \frac{1}{4}\begin{pmatrix} n & K_0 n^{1-\widehat{\alpha}_{n_0}} \\ K_0 n^{1-\widehat{\alpha}_{n_0}} & K_0^2 a_n(\widehat{\alpha}_{n_0}) \end{pmatrix}.$$

Therefore,

$$\mathbf{\Psi}_n\left[\sum_{k=1}^{n} \mathbf{f}'_k(\theta_0,\widehat{\nu})f'^T_k(\theta_0,\widehat{\nu})\right]^{-1}\sum_{k=1}^{n} \eta_k \mathbf{f}'_k(\theta_0) = \Phi_n^{-1}\sum_{k=1}^{n} \eta_k \mathbf{f}'_k(\theta_0,\widehat{\nu}).$$

Then, using the CLT for martingale arrays ([4, 29]), we have

$$\lim_n \mathbf{\Psi}_n(\widehat{\theta}_{n_0} - \theta_0) \stackrel{\mathcal{D}}{=} \mathcal{N}(0, (K_0/2)\mathbf{I}) \iff \lim_n \Phi_n(\widehat{\theta}_{n_0} - \theta_0) \stackrel{\mathcal{D}}{=} \mathcal{N}(0, (K_0/2)\mathbf{I}).$$

**9. Extension to estimating functions.** To simplify notation, we assume that $q = 0$. Let $\mathbf{Q}_n(\theta; \{Z_k, \mathbf{a}_k(\theta)\})$ be defined as in (1.3), where $\mathbf{a}_k(\theta) := \mathbf{g}'_k(\theta)[b_k(\theta)]^{-1}$, $b_k(\theta)$ being $\mathcal{F}_{k-1}$-measurable. Let $\widehat{\theta}_n$ solve $\mathbf{Q}_n(\theta; \{Z_k, \mathbf{a}_k(\theta)\}) = 0$. The CLSE solution of (1.1) corresponds to $b_k(\theta) = \lambda_k$ and the QLE to $b_k(\theta) = \sigma_k^{-2}(\theta)$.

Following [12], let

$$(9.1) \qquad S_n(\theta) := \left[ \int Q_{n;i}(u; \{Z_k, \mathbf{a}_k(\theta)\}) \, du_i \right](\theta), \qquad i = 1, \ldots, n,$$

when this quantity exists. Since $\int Q_{n;i}(u; \{Z_k, \mathbf{a}_k(\theta)\}) \, du_i](\theta) = -\sum_{k=1}^n \int (Z_k - g_k(\theta)) g'_{k;i}(\theta)[b_k(\theta)]^{-1}$, it follows that when $b_k(\theta)$ is a function of $g_k(\theta)$, $S_n(\theta)$ exists (but does not necessarily have an explicit form).

Let $e_k(\theta) := Z_k - g_k(\theta)$, $d_k^g(\theta) := g_k(\theta_0) - g_k(\theta)$, $D_{a,n}(\theta) := \sum_{k=1}^n d_{a,k}(\theta)$, where, for $i = 1, \ldots, p$,

$$d_{a,k}(\theta) := \left[ \int d_k^g(u) g'_{k;i}(u) b_k^{-1}(u) \, du_i \right](\theta_0) - \left[ \int d_k^g(u) g'_{k;i}(u) b_k^{-1}(u) \, du_i \right](\theta)$$

and let $d_{b,k}(\theta) = [\int g'_{k;i}(u) b_k^{-1}(u) \, du_i](\theta_0) - [\int g'_{k;i}(u) b_k^{-1}(u) \, du_i](\theta)$, for $i = 1, \ldots, p$, when these quantities exist.

PROPOSITION 9.1. *Let us assume that $S_n(\theta)$, defined as in (9.1), exists and let $\Omega_\infty \subset \Omega$ be defined by*

$$\Omega_\infty = \{\mathrm{LIP}_\theta(\{d_{a,k}(\theta)\}) \cap \mathrm{LIP}_\theta(\{d_{b,k}(\theta)\}) \cap \mathrm{SI}_\theta(\{D_{a,n}(\theta)\})$$
$$\cap \mathrm{VAR}_\theta(\{\sigma_k^2(\theta_0, \nu_0), d_{b,k}(\theta), D_{a,k}(\theta)\})\}.$$

*Let us assume that $P(\Omega_\infty) > 0$. Then $\lim_n \widehat{\theta}_n \stackrel{\text{a.s.}}{=} \theta_0$ on $\Omega_\infty$.*

The proof is the same as in Section 3, where $S_n(\theta) - S_n(\theta_0) = \sum_{k=1}^n d_{a,k} + \sum_{k=1}^n e_k d_{b,k}(\theta)$. In the particular case $b_k(\theta) = \lambda_k$, we have $d_{a,k} = 2^{-1}[d_k^g(\theta)]^2 \lambda_k^{-1}$, $d_{b,k}(\theta) = d_k^g(\theta) \lambda_k^{-1}$ and $D_{a,n}(\theta) = D_n(\theta)$ up to additive constants, where $D_n(\theta)$ is given by (1.8). $\Omega_\infty$ is then reduced to the subset defined in Proposition 3.1.

In the general case where $b_k(\theta)$ depends on $\theta$, since conditions given in Proposition 9.1 may be difficult or even impossible to verify in practice, we may derive sufficient conditions using Wu's lemma applied to a Taylor series expansion of $S_n(\theta) - S_n(\theta_0)$. Writing $e_k(\theta_0) =: e_k$, $\theta_a = \theta_0 + a(\theta - \theta_0)$, $a \in (0, 1)$, we have

$$c_k(\theta) := (\theta - \theta_0)^T \mathbf{g}'_k(\theta_0) b_k^{-1}(\theta_0) + c_{a,k}(\theta),$$
$$c_{a,k}(\theta) := (\theta - \theta_0)^T [\mathbf{g}''_k(\theta_a) b_k^{-1}(\theta_a) - \mathbf{g}'_k(\theta_a)[\mathbf{b}'_k(\theta_a)]^T b_k^{-2}(\theta_a)](\theta - \theta_0),$$
$$d_{a,k}(\theta) := (\theta - \theta_0)^T \mathbf{g}'_k(\theta_a),$$

$$D_{a,n}(\theta) := \sum_{k=1}^{n} [d_{a,k}(\theta)]^2 = (\theta - \theta_0)^T \sum_{k=1}^{n} \mathbf{g}'_k(\theta_a)[\mathbf{g}'_k(\theta_a)]^T(\theta - \theta_0).$$

PROPOSITION 9.2. *Let us assume that $S_n(\theta)$ defined as in (9.1) exists and that $g_n(\theta)$ has second derivatives in $\theta$. Let $\Omega_\infty \subset \Omega$ be defined by*

$$\Omega_\infty = \mathrm{LIP}_\theta(\{c_{a,k}(\theta)\}) \cap \bigcap_i \mathrm{LIP}_\theta(\{g'_{k;i}(\theta)\})$$

$$\cap \mathrm{SI}_\theta(\{D_{a,n}(\theta)\}) \cap \mathrm{VAR}_\theta(\{\sigma_k^2(\theta_0, \nu_0), c_k(\theta), D_{a,k}(\theta)\})$$

$$\cap \left\{ \overline{\lim_n} \sup_\theta \left[ \sum_{k=1}^{n} d_k^g(\theta) c_{a,k}(\theta) \right] [D_{a,n}(\theta)]^{-1} \stackrel{\mathrm{a.s.}}{=} 0 \right\}.$$

*Let us assume that $P(\Omega_\infty) > 0$. Then $\lim_n \widehat{\theta}_n \stackrel{\mathrm{a.s.}}{=} \theta_0$ on $\Omega_\infty$.*

In the particular class $b_k(\theta) = \lambda_k$, we have $\mathrm{VAR}_\theta(\cdot) = \mathrm{VAR}_\theta(\{\sigma_k^2, (\theta - \theta_0)^T \mathbf{f}'_k(\theta), D_{a,k}(\theta)\}) = \mathrm{VAR}_\theta(\{\sigma_k^2, d_k(\theta), D_k(\theta)\})$, where $d_k(\theta), D_k(\theta)$ are given by (1.8). Consequently, $\Omega_\infty$ (Proposition 9.2) $\subset \Omega_\infty$ (Proposition 3.1). In the linear case with $b_k(\theta) = \lambda_k$, we have that $\Omega_\infty$ is reduced to (1.7) in both propositions.

PROOF OF PROPOSITION 9.2. Let us write $S_n(\theta) - S_n(\theta_0) = (\theta - \theta_0)^T \times \mathbf{S}'_n(\theta_0) + (\theta - \theta_0)^T \mathbf{S}''_n(\theta_a)(\theta - \theta_0)$. Then

$$S_n(\theta) - S_n(\theta_0)$$

$$= D_{a,n}(\theta) - (\theta - \theta_0)^T \sum_{k=1}^{n} e_k \mathbf{g}'_k(\theta_0) b_k^{-1}(\theta_0)$$

$$- (\theta - \theta_0)^T \sum_{k=1}^{n} e_k [\mathbf{g}''_k(\theta_a) b_k^{-1}(\theta_a) - \mathbf{g}'_k(\theta_a)[b'_k(\theta_a)]^T b_k^{-2}(\theta_a)](\theta - \theta_0)$$

$$- (\theta - \theta_0)^T \sum_{k=1}^{n} d_k^g(\theta_a) [\mathbf{g}''_k(\theta_a) b_k^{-1}(\theta_a) - \mathbf{g}'_k(\theta_a)[b'_k(\theta_a)]^T b_k^{-2}(\theta_a)]$$

$$\times (\theta - \theta_0)$$

and then the proof is as in Section 3. □

We will not describe here the asymptotic distribution for these estimators since the methodology is the same as for the CLSE.

# APPENDIX

PROOF OF PROPOSITION 5.1. We work on $\mathrm{VAR}_\theta(\{\sigma_k^2, d_k(\theta), D_{*k}(\theta)\})$. Let us assume that $\widetilde{\Theta}$ is a finite set (i.e., $\widetilde{\Theta} = \{\theta_i\}_{i \leq I}$). Therefore, there exist some random integers $i_n \leq I$ and $j_n \leq I$ such that $\sup_{\theta \in \widetilde{\Theta}} |L_n(\theta)| [D_{*n}(\theta)]^{-1} \stackrel{\mathrm{a.s.}}{=} |L_n(\theta_{i_n})| [D_{*n}(\theta_{i_n})]^{-1}$, $\inf_{\theta \in \widetilde{\Theta}} D_{*n}(\theta) \stackrel{\mathrm{a.s.}}{=} D_{*n}(\theta_{j_n})$. Moreover, thanks to the SLLNM (strong law of large numbers for martingales, Theorem 2.18, [9]), we have

$$\forall \theta \in \widetilde{\Theta} \qquad \lim_n |L_n(\theta)| [D_{*n}(\theta)]^{-1} \stackrel{\mathrm{a.s.}}{=} 0 \qquad \text{on } \left\{ \lim_n D_{*n}(\theta) \stackrel{\mathrm{a.s.}}{=} \infty \right\},$$

which implies that $\lim_n |L_n(\theta_{i_n})| [D_{*n}(\theta_{i_n})]^{-1} \stackrel{\mathrm{a.s.}}{=} 0$ on $\bigcap_{i \leq I} \{\lim_n D_{*n}(\theta_i) \stackrel{\mathrm{a.s.}}{=} \infty\}$, which contains $\{\lim_n D_{*n}(\theta_{j_n}) \stackrel{\mathrm{a.s.}}{=} \infty\}$.

Let us now assume the general case $\widetilde{\Theta} \subset \mathbb{R}^p$. The general idea for proving the result is then to extend the proof concerning the case where $\widetilde{\Theta}$ is finite: for each $n$, we will use some random $\mathcal{F}_{n-1}$-measurable discretization $\widetilde{\Theta}_n$ of $\widetilde{\Theta}$ that becomes finer and finer as $n \to \infty$, together with the fact that, for each $k \leq n$, there will exist a point of $\widetilde{\Theta}_k$ which will get closer and closer to $\theta$ as $k$ increases. So, writing

$$(\mathrm{A.1}) \quad \overline{\lim}_n \sup_\theta \frac{|L_n(\theta)|}{D_{*n}(\theta)} \leq \overline{\lim}_n \sup_\theta \frac{|L_n(\theta, \mathcal{G}_n(\theta))|}{D_{*n}(\theta)} + \overline{\lim}_n \sup_\theta \frac{|L_n(\mathcal{G}_n(\theta))|}{D_{*n}(\theta)},$$

where $\mathcal{G}_n(\theta) := \{\theta_k(\theta)\}_{k \leq n}$, $\theta_k(\theta) \in \widetilde{\Theta}_k$ being one vertex of $\widetilde{\Theta}_k$ among those closest to $\theta$, and

$$L_n(\theta, \mathcal{G}_n(\theta)) := \sum_{k=1}^n \eta_k (d_k(\theta) - d_k(\theta_k(\theta))), \qquad L_n(\mathcal{G}_n(\theta)) := \sum_{k=1}^n \eta_k d_k(\theta_k(\theta)),$$

it will be sufficient to prove that each of the terms of the sum of the right-hand side of (A.1) is null. For the second term, we will use the fact that the set $\{\{\theta_k(\theta)\}_{k \leq n}\}_\theta$ is finite, together with the usual SLLNM (strong law of large numbers for martingales), $\mathrm{LIP}_\theta(\{d_{*k}(\theta)\})$ and the sufficiently rapid convergence to 0 of the mesh size. For the first term, we define

$$(\mathrm{A.2}) \qquad U_{m,n}(\theta, \mathcal{G}_n(\theta)) := \sum_{k=m}^n \eta_k (d_k(\theta) - d_k(\theta_k(\theta))) [D_{*k}(\theta)]^{-1}$$

and we will use a property of submartingales (Theorem 2.1 in [9]; see also Theorem A.1) that leads to

$$(\mathrm{A.3}) \quad \begin{aligned} \lambda P\Big( \max_{n : m \leq n \leq m'} \sup_\theta |U_{m,n}(\theta, \mathcal{G}_n(\theta))| > \lambda \Big) \\ \leq E\Big( \sup_\theta |U_{m,m'}(\theta, \mathcal{G}_{m'}(\theta))| \Big). \end{aligned}$$

Then, using $\mathrm{LIP}_\theta(\{d_k(\theta)\})$ and the sufficiently rapid convergence of the mesh size to 0, we will prove that $\lim_m \lim_{m'} E(\sup_\theta |U_{m,m'}(\theta, \mathcal{G}_{m'}(\theta))|) = 0$, from which we will deduce that $\lim_m \sup_{n \geq m} \sup_\theta |U_{m,n}(\theta, \mathcal{G}_n(\theta))| \stackrel{P}{=} 0$, thanks to (A.3), and then $\lim_m \overline{\lim}_{n \geq m} \sup_\theta |U_{m,n}(\theta, \mathcal{G}_n(\theta))| \stackrel{\mathrm{a.s.}}{=} 0$. Finally, the result will follow from the relationship, due to Lemma A.3,

$$|L_n(\theta, \mathcal{G}_n(\theta))|[D_{*n}(\theta)]^{-1} = \left|\sum_{k=1}^n U_{k,n}(\theta, \mathcal{G}_n(\theta))[d_{*k}(\theta)]^2\right|[D_{*n}(\theta)]^{-1},$$

together with a generalized Toeplitz lemma applied to the $\sup_\theta$ of this quantity.

We now provide details of the proof. We define, for each $k$, a discretization of $\mathbb{R}^p$ by a random grid $G_k$ with fixed directions, a fixed origin and a random mesh size $\epsilon_k$, $\mathcal{F}_{k-1}$-measurable and converging a.s. to 0 sufficiently rapidly as $k \to \infty$ according to the assumptions $A_1^\epsilon$ and $A_2^\epsilon$ defined later in the proof. Let $\{\theta_{k,i}\}_i =: \widetilde{\Theta}_k$ be the vertices of $G_k \cap \widetilde{\Theta}$ and let, for $\theta \in \widetilde{\Theta}$, $\theta_k(\theta) \in \widetilde{\Theta}_k$ be one of the elements of $\widetilde{\Theta}_k$ closest to $\theta$, that is, $\|\theta_k(\theta) - \theta\| \leq c\epsilon_k$, where $c$ is any constant satisfying $c \geq \sqrt{p}/2$ where $\|\cdot\|$ is the Euclidean norm.

Let us first consider the second term of the right-hand side of (A.1). Defining $D_n(\mathcal{G}_n(\theta)) := \sum_{k=1}^n [d_k(\theta_k(\theta))]^2$, $D_{*n}(\mathcal{G}_n(\theta)) := \sum_{k=1}^n [d_{*k}(\theta_k(\theta))]^2$, we have

(A.4)
$$\sup_\theta \frac{|L_n(\mathcal{G}_n(\theta))|}{D_{*n}(\theta)}$$
$$\leq \sup_\theta \frac{|L_n(\mathcal{G}_n(\theta))|}{D_{*n}(\mathcal{G}_n(\theta))}\left[\sup_\theta \frac{|D_{*n}(\mathcal{G}_n(\theta)) - D_{*n}(\theta)|}{D_{*n}(\theta)} + 1\right].$$

Since, for $k \leq n$, $\theta_k(\theta)$ is $\mathcal{F}_{k-1}$-measurable and $\{\mathcal{G}_n(\theta)\}_\theta$ is a finite set, we get, as in the finite $\widetilde{\Theta}$ case,

(A.5) $\quad \limsup_n \sup_\theta \frac{|L_n(\mathcal{G}_n(\theta))|}{D_{*n}(\mathcal{G}_n(\theta))} \stackrel{\mathrm{a.s.}}{=} 0 \quad$ on $\left\{\liminf_n \inf_\theta D_{*n}(\mathcal{G}_n(\theta)) \stackrel{\mathrm{a.s.}}{=} \infty\right\}.$

Moreover, thanks to $\mathrm{LIP}_\theta(\{d_{*k}(\theta)\})$, there exists $g_{*k}$, $\mathcal{F}_{k-1}$-measurable, and $h_*(\cdot)$ such that $|d_{*k}(\theta) - d_{*k}(\theta_k(\theta))| \leq h_*(\|\theta_k(\theta) - \theta_k\|)g_{*k}$, implying that

(A.6) $\quad \sup_\theta \frac{|D_{*n}(\mathcal{G}_n(\theta)) - D_{*n}(\theta)|}{D_{*n}(\theta)} \leq \frac{\sum_{k=1}^n \varepsilon_{*k} g_{*k} u_k(\varepsilon_{*k})}{\inf_\theta D_{*n}(\theta)},$

where $\varepsilon_{*k} := h_*(c\epsilon_k)$ and $u(\varepsilon_{*k}) = 2\max_i |d_{*k}(\theta_{k,i})| + \varepsilon_{*k} g_{*k}$.

Let $0 < a < 1$ and let us choose $\{\epsilon_k\}$ such that $A_1^\epsilon : \varepsilon_{*k} g_{*k} u_k(\varepsilon_{*k}) \stackrel{\mathrm{a.s.}}{<} a^k$ for all $k$. Therefore, the limit in $n$ of (A.6) is a.s. finite. Then, also using (A.4) and (A.5), we get

$$\limsup_n \sup_\theta |L_n(\mathcal{G}_n(\theta))|[D_{*n}(\theta)]^{-1} \stackrel{\mathrm{a.s.}}{=} 0 \quad \text{on } \left\{\liminf_n \inf_\theta D_{*n}(\mathcal{G}_n(\theta)) \stackrel{\mathrm{a.s.}}{=} \infty\right\}.$$

Moreover, $\lim_n \inf_\theta D_{*n}(\mathcal{G}_n(\theta)) \stackrel{\text{a.s.}}{=} \infty$ is equivalent to $\lim_n \inf_\theta D_{*n}(\theta) \stackrel{\text{a.s.}}{=} \infty$ under $A_1^\epsilon$ because

$$\liminf_n \inf_\theta D_{*n}(\theta) \le \lim_n \inf_\theta |D_{*n}(\theta) - D_{*n}(\mathcal{G}_n(\theta))| + \liminf_n \inf_\theta D_{*n}(\mathcal{G}_n(\theta))$$

$$\le \sum_{k=1}^\infty \varepsilon_{*k} g_{*k} u_k(\varepsilon_{*k}) + \liminf_n \inf_\theta D_{*n}(\mathcal{G}_n(\theta))$$

and $\lim_n \inf_\theta D_{*n}(\mathcal{G}_n(\theta)) \le \sum_{k=1}^\infty \varepsilon_{*k} g_{*k} k u_k(\varepsilon_{*k}) + \lim_n \inf_\theta D_{*n}(\theta)$.

Next, we show that $\lim_n \sup_\theta |L_n(\theta, \mathcal{G}_n(\theta))|[D_{*n}(\theta)]^{-1} \stackrel{\text{a.s.}}{=} 0$. Let us write

$$(d_k(\theta) - d_k(\theta_k(\theta)))[D_{*k}(\theta)]^{-1} =: \widetilde{d}_k(\theta, \theta_k(\theta)).$$

Since $\{\eta_k\}$ is a martingale difference sequence and $d_k(\theta)$, $\theta_k(\theta)$ and $D_k(\theta)$ are $\mathcal{F}_{k-1}$-measurable, it follows that $\{U_{m,n}(\theta, \mathcal{G}_n(\theta))\}_n$, defined as in (A.2), is a martingale:

$$U_{m,n-1}(\theta, \mathcal{G}_{n-1}(\theta)) = E(U_{m,n}(\theta, \mathcal{G}_n(\theta))|\mathcal{F}_{n-1}).$$

According to Jensen's inequality, this implies that $\{\sup_\theta |U_{m,n}(\theta, \mathcal{G}_n(\theta))|\}_n$ is a submartingale:

$$\sup_\theta |U_{m,n-1}(\theta, \mathcal{G}_{n-1}(\theta))| \le E\Big(\sup_\theta |U_{m,n}(\theta, \mathcal{G}_n(\theta))| \Big| \mathcal{F}_{n-1}\Big).$$

Therefore, using Theorem 2.1 ([9]; see also Theorem A.1) and denoting by $V_{m,n}$ the quantity $\sup_\theta |U_{m,n}(\theta, \mathcal{G}_n(\theta))|$, we get, for any $\lambda > 0$,

$$(\text{A.7}) \qquad \lambda P\Big(\max_{n:\, m \le n \le m'} V_{m,n} > \lambda\Big) \le E(V_{m,m'}).$$

Now, using $V_{m,m'} := \sup_\theta |\sum_{k=m}^{m'} \eta_k \widetilde{d}_k(\theta, \theta_k(\theta))|$ and $\text{LIP}_\theta(\{d_k(\theta)\})$, which implies that $|\widetilde{d}_k(\theta, \theta_k(\theta))| \le \varepsilon_k g_k [D_{*k}(\theta)]^{-1}$, where $g_k$ is $\mathcal{F}_{k-1}$-measurable and $\varepsilon_k = h(c\epsilon_k)$, we get

$$E(V_{m,m'}) \le E\bigg(\sum_{k=m}^{m'} E\Big(|\eta_k|\varepsilon_k g_k \Big[\inf_\theta D_{*k}(\theta)\Big]^{-1} \Big| \mathcal{F}_{k-1}\Big)\bigg).$$

Using Hölder's inequality and the $\mathcal{F}_{k-1}$-measurability of $\varepsilon_k$, $g_k$, $[\inf_\theta D_{*k}(\theta)]^{-1}$, we get

$$E(V_{m,m'}) \le E\bigg(\sum_{k=m}^{m'} \sigma_k \varepsilon_k g_k \Big[\inf_\theta D_{*k}(\theta)\Big]^{-1}\bigg).$$

Let us choose $\{\epsilon_k\}$ such that it satisfies $A_2^\epsilon : \varepsilon_k \sigma_k g_k [\inf_\theta D_{*k}(\theta)]^{-1} \stackrel{\text{a.s.}}{\le} a^k$, for all $k$, in addition to $A_1^\epsilon$. Then $E(V_{m,m'}) \le \sum_{k=m}^\infty a^k < \infty$, implying, according to (A.7), that

$$(\text{A.8}) \qquad \lim_{m'} \lambda P\Big(\max_{n:\, m \le n \le m'} V_{m,n} > \lambda\Big) \le \sum_{k=m}^\infty a^k.$$

Consequently, since $\sum_{k=m}^{\infty} a^k < \infty$, (A.8) implies that

$$\text{(A.9)} \qquad \lim_m \lim_{m'} P\Big(\max_{n:\, m \leq n \leq m'} V_{m,n} > \lambda\Big) = 0.$$

Moreover, since $\{\max_{n:\, m \leq n \leq m'} V_{m,n} > \lambda\}_{m'}$ is an increasing sequence of events, it follows that $P(\sup_{n:\, m \leq n} V_{m,n} > \lambda) = \lim_{m'} P(\max_{m \leq n \leq m'} V_{m,n} > \lambda)$. Therefore, (A.9) becomes $\lim_m P(\sup_{n \geq m} V_{m,n} > \lambda) = 0$ for all $\lambda > 0$, which means that $\lim_m \sup_{n \geq m} V_{m,n} \stackrel{P}{=} 0$. Therefore, there is a subsequence $\{\sup_{n \geq m_i} V_{m_i,n}\}_{m_i}$ that converges a.s. to 0 as $m_i \to \infty$. However, for $m > m_i$ with $m \leq n$, $U_{m,n}(\theta, \mathcal{G}_n(\theta)) = U_{m_i,n}(\theta, \mathcal{G}_n(\theta)) - U_{m_i,m-1}(\theta, \mathcal{G}_{m-1}(\theta))$, which implies that $\sup_{n \geq m} V_{m,n} \leq 2\sup_{n \geq m_i} V_{m_i,n}$ and, consequently,

$$\text{(A.10)} \quad \lim_m \sup_{n \geq m} V_{m,n} \stackrel{\text{a.s.}}{=} 0 \qquad \text{which implies that } \varlimsup_m \varlimsup_{n \geq m} V_{m,n} \stackrel{\text{a.s.}}{=} 0.$$

It then remains to deduce from (A.10) that $\lim_n \sup_\theta L_n(\theta, \mathcal{G}_n(\theta))[D_{*n}(\theta)]^{-1} \stackrel{\text{a.s.}}{=} 0$. Let us write $S_k$ for $\sum_{l=1}^{k-1} \eta_l \widetilde{d}_l(\theta, \theta_l(\theta)) := U_{1,k-1}(\theta, \mathcal{G}_{k-1}(\theta))$. Then

$$L_n(\theta, \mathcal{G}_n(\theta)) = \sum_{k=1}^n (S_{k+1} - S_k) D_{*k}(\theta).$$

Using Lemma A.3 and $D_{*k}(\theta) - D_{*k-1}(\theta) = [d_{*k}(\theta)]^2$, we get

$$L_n(\theta, \mathcal{G}_n(\theta)) = \sum_{k=1}^n (S_{n+1} - S_k)[d_{*k}(\theta)]^2 = \sum_{k=1}^n U_{k,n}(\theta, \mathcal{G}_n(\theta))[d_{*k}(\theta)]^2,$$

implying that

$$\text{(A.11)} \quad \varlimsup_n \sup_\theta \frac{|L_n(\theta, \mathcal{G}_n(\theta))|}{D_{*n}(\theta)} \leq \varlimsup_N \varlimsup_{n \geq N} \sup_{k < N} V_{k,n} \varlimsup_N \varlimsup_{n \geq N} \sup_\theta \frac{D_{*N}(\theta)}{D_{*n}(\theta)}$$
$$+ \varlimsup_N \varlimsup_{n \geq N} \sup_{N < k \leq n} V_{k,n}.$$

Now, using, in the first term of the right-hand of side of (A.11), $U_{k,n}(\theta, \mathcal{G}_n(\theta)) = U_{k,N-1}(\theta, \mathcal{G}_{N-1}(\theta)) + U_{N,n}(\theta, \mathcal{G}_n(\theta))$, and, in the second term, $U_{k,n}(\theta, \mathcal{G}_n(\theta)) = U_{N,n}(\theta, \mathcal{G}_n(\theta)) - U_{N,k-1}(\theta, \mathcal{G}_{k-1}(\theta))$, thanks to (A.10), we get (5.1). □

LEMMA A.1 (Wu's lemma [34]). *If for all $\delta > 0$ sufficiently small, $\varliminf_n \inf_{\theta \in B_\delta^c}(S_n(\theta) - S_n(\theta_0)) \stackrel{\text{a.s.}}{>} {}^{(P)} 0$, then $\lim_n \widehat{\theta}_n \stackrel{\text{a.s.}}{=} {}^{(P)} \theta_0$.*

THEOREM A.1 (Theorem 2.1 in [9]). *If $\{S_i, \mathcal{F}_i, 1 \leq i \leq n\}$ is a submartingale, then, for each real $\lambda$, we have $\lambda P(\max_{i \leq n} S_i > \lambda) \leq E(S_n 1_{\{\max_{i \leq n} S_i > \lambda\}})$.*

LEMMA A.2 [16]. *Let $a_k \geq 0$ for all $k$, with $a_1 > 0$, and $S_n = \sum_{k=1}^n a_k$ with $\lim_n S_n \leq \infty$. Then $\sum_{k=1}^\infty a_k S_k^{-2} \leq 2a_1^{-1} - \lim_n S_n^{-1}$.*

LEMMA A.3. *Let $\{S_k\}$ with $S_1 = 0$, $D_k = \sum_{l=1}^k d_l^2$. Then $\sum_{k=1}^n (S_{k+1} - S_k)D_k = \sum_{k=1}^n (S_{n+1} - S_k)d_k^2$.*

PROOF. $\sum_{k=1}^n (S_{k+1} - S_k)D_k = S_{n+1}D_n + \sum_{k=1}^{n-1} S_{k+1}D_k - \sum_{k=1}^n S_k D_k = S_{n+1}D_n + \sum_{k=1}^n S_k(D_{k-1} - D_k)$. Then use $D_k - D_{k-1} = d_k^2$. $\square$

**Acknowledgments.** The author is grateful to V. Vatutin and A. Zubkov from the Steklov Institute (Moscow) for their fruitful comments.

UR341  
National Agronomical  
　Research Institute (INRA)  
F-78352 Jouy-en-Josas  
France  
E-mail: christine.jacob@jouy.inra.fr